\documentclass[a4paper,twoside,11pt]{article}

\usepackage{amsfonts}
\usepackage{mathrsfs}
\usepackage{amsmath}
\usepackage{amsmath,graphics}
\usepackage{amssymb}
\usepackage{indentfirst,latexsym,bm,amsmath,amssymb,amsthm,graphicx,fancyhdr,float}
\usepackage[all]{xy}
\usepackage{amsthm}
\usepackage{amscd}
\usepackage{tikz}
\usetikzlibrary{cd}
\usepackage[hyperfootnotes=false,colorlinks,linkcolor=black,citecolor=black,urlcolor=black]{hyperref}
\usepackage[center]{titlesec}
\usepackage{booktabs}
\usepackage{xcolor}
\usepackage{adjustbox}
\newcommand{\papertablestyle}{%
  \footnotesize\renewcommand{\arraystretch}{1.48}%
  \setlength{\tabcolsep}{5pt}%
}
\titleformat{\section}{\centering\Large\bfseries}{\thesection.}{1em}{}
\titleformat{\subsection}{\centering\large\bfseries}{\thesubsection.}{1em}{}

\begin{document}
\theoremstyle{plain}
\newtheorem{thm}{Theorem}
\newtheorem{theorem}{Theorem}[section]
\newtheorem{pro}[theorem]{Proposition}
\newtheorem{lemma}[theorem]{Lemma}
\newtheorem{corollary}[theorem]{Corollary}
\newtheorem{proposition}[theorem]{Proposition}
\newtheorem{conjecture}[theorem]{Conjecture}
\newtheorem{problem}[theorem]{Problem}
\theoremstyle{definition}
\newtheorem{remark}[theorem]{Remark}
\newtheorem{definition}[theorem]{Definition}
\newtheorem{claim}{Claim}[theorem]
\newtheorem{exa}[theorem]{Example}
\newtheorem{news}[theorem]{}

\numberwithin{equation}{section}
\def\wt{\widetilde}
\def\ol{\overline}
\def\ra{\rightarrow}
\def\lra{\longrightarrow}
\def\div{\text{\rm{div\,}}}
\def\Div{\text{\rm{Div\,}}}
\def\Aut{\text{\rm{Aut\,}}}
\def\({$($}
\def\){$)$}
\def\chit{\chi_{\rm top}}
\def\bbp{\mathbb P}
\def\calh{\mathcal H}
\def\cala{\mathcal A}
\def\calt{\mathcal T}
\def\calm{\mathcal M}
\def\caln{\mathcal N}
\def\cals{\mathcal S}
\def\calf{\mathcal F}
\def\Pic{\text{{\rm Pic\,}}}
\def\rank{\text{{\rm rank\,}}}

\font\eightrm=cmr10 at 8pt \font\bigrm=cmr10 scaled1440
\font\midrm=cmr10 scaled1200 \font\ninerm=cmr10 at 9pt
\font\ninebf=cmbx10 at 9pt \font\bbf=msbm10 \font\bigbf=cmbx10
scaled1440 \font\Bigbf=cmbx10 scaled2073 \font\midbf=cmbx10
scaled1200 \font\eightit=cmti8 \font\nineit=cmti9
\renewcommand{\labelenumi}{$($\arabic{enumi}$)$}
\renewcommand{\thethm}{\Alph{thm}}
\renewcommand{\theequation}{\arabic{section}-\arabic{equation}}
\renewcommand{\labelenumii}{(\Alph{enumi})}
\pagestyle{fancy}
\setlength{\headsep}{0.4cm}
\fancyhf{}
\fancyhead[CE]{\small Xin L\"u and Wan-Yuan Xu}
\fancyhead[CO]{\small Negative Segre number and big cotangent bundle}
\fancyfoot[C]{\small\thepage}
\renewcommand{\headrulewidth}{0.4pt}
\renewcommand{\footrulewidth}{0pt}
\fancypagestyle{plain}{%
  \fancyhf{}
  \fancyfoot[C]{\small\thepage}
  \renewcommand{\headrulewidth}{0pt}
  \renewcommand{\footrulewidth}{0pt}
}

\title{ Orbifold Chern slopes of surfaces with negative second Segre number and big cotangent bundle}

\author{Xin L\"u, Wan-Yuan Xu}
\date{}
\renewcommand{\thefootnote}{}

\maketitle

\begin{abstract}
In this paper, we construct smooth minimal projective surfaces of general type with negative second Segre number and big cotangent bundle, and prove that every rational number in $(1,3]$ is realized as the orbifold Chern slope of a Du Val contraction of such a surface. The upper endpoint is attained by a canonical model with eighteen $A_8$ singularities and three $A_2$ singularities, which is a compact complex two-ball quotient by the equality case of the orbifold Miyaoka–Yau inequality. We also construct an explicit sequence of canonical-model orbifold Chern slopes converging to three from below. Our constructions are based on ramified base changes of semistable fibrations.
\end{abstract}

\section{Introduction}\label{outline}
We work over the field of complex numbers $\mathbb C$ throughout this paper.
Let $Y$ be a smooth minimal projective surface of general type, and put
\[
 s_2(Y)=c_1^2(Y)-c_2(Y)=K_Y^2-\chi_{\rm top}(Y).
\]
Bogomolov~\cite{Bo} proved that $s_2(Y)>0$ implies that the cotangent
bundle $\Omega_Y^1$ is big. A surface of general type with big
cotangent bundle contains only finitely many rational or elliptic
curves. It is natural to ask whether the cotangent bundle can still
be big when $s_2(Y)<0$.

The question is also related to the existence of symmetric differentials
and to the hyperbolicity of projective varieties; see
\cite{Br17,DMR,DR,RR,R21}. For example, an ample cotangent bundle
implies Kobayashi and algebraic hyperbolicity
\cite{De04,De05,Kob98}. On the other hand, a smooth hypersurface in
$\mathbb P^3$ has no nonzero symmetric differentials of positive
order~\cite{Br71}. For a smooth hypersurface of degree $d\geq5$,
one has $s_2=10d-4d^2<0$. Related vanishing results for weighted
hypersurfaces were obtained in~\cite{Xu26}.

Bogomolov and De Oliveira~\cite{BD} showed that resolutions of
nodal hypersurfaces may acquire symmetric differentials even though
the corresponding smooth hypersurfaces have none.
Roulleau and Rousseau~\cite{RR14} proved that
\[
 s_2(Y)+s_2(\mathcal Z)>0
\]
is sufficient for the cotangent bundle of the minimal resolution
$Y$ of a canonical surface $Z$ to be big, where $\mathcal Z$ is
the natural orbifold of $Z$. Their construction gives degree-$13$
hypersurfaces whose minimal resolutions have negative second Segre
number and big cotangent bundle. More recently, Asega, De Oliveira
and Weiss~\cite{AOW0,AOW} established the QS-bigness criterion,
where QS abbreviates \emph{quotient singularities}.
This criterion also applies to quotients that are not canonical
models and yields examples from degree-$8$ hypersurfaces; see
\cite[Theorem~1.5]{AOW0}.

We are interested in the geography of such surfaces. Let $Z$ be a
projective normal surface with only Du Val singularities, let
$Y\to Z$ be its minimal resolution, and denote by $\mathcal Z$ its
natural orbifold. When $c_2(\mathcal Z)>0$, we call
\[
 \mathrm{slope}(\mathcal Z)=\frac{c_1^2(\mathcal Z)}{c_2(\mathcal Z)}
\]
the orbifold Chern slope. If $Y$ is minimal and of general type,
then $K_Z$ is nef and big, and the orbifold Bogomolov--Miyaoka--Yau
inequality gives $0<\mathrm{slope}(\mathcal Z)\leq3$; see
\cite{Miy84} and \cite[Theorem~1.1]{GKPT19}. In the positive orbifold
Segre region, namely
\[
 s_2(\mathcal Z):=c_1^2(\mathcal Z)-c_2(\mathcal Z)>0,
\]
the slope therefore belongs to $(1,3]\cap\,\mathbb Q$.
The lower bound is equivalent to this positivity condition.
For the $A_n$ quotients satisfying the QS inequality used below,
the condition follows directly from the local estimates; see
Proposition~\ref{prop:orbifold-slope-range}.

In~\cite[Section~5]{RR14}, Roulleau and Rousseau considered how
close this slope can be to three when $s_2(Y)<0$ and $\Omega_Y^1$
is big. Our main result realizes every rational slope in the
positive orbifold Segre range.

\begin{theorem}\label{thm:full-rational-interval}
For every rational number $r\in(1,3]$, there exists a projective
normal surface $Z$ with only Du Val singularities whose minimal
resolution $Y$ is a smooth minimal surface of general type such that
\[
 s_2(Y)<0,\qquad \Omega_Y^1\text{ is big},\qquad
 \frac{c_1^2(\mathcal Z)}{c_2(\mathcal Z)}=r.
\]
Here $\mathcal Z$ is the natural orbifold of $Z$.
\end{theorem}

Consequently, among Du Val quotients with $s_2(\mathcal Z)>0$
whose minimal resolutions satisfy the conditions of the theorem,
the set of orbifold Chern slopes is exactly $(1,3]\cap\mathbb Q$.
The positivity condition concerns the chosen quotient; it is not
implied by $s_2(Y)<0$ and the bigness of $\Omega_Y^1$ alone.

The upper endpoint is attained by a canonical model. More precisely,
we obtain the following surface from the genus-two fibration of
Tu~\cite[Theorem~1.1(2) and Section~3]{Tu07}.

\begin{theorem}\label{thm:sharp-three}
A connected cyclic base change of degree nine of the fibration in
Subsection~\ref{subsec:tu-five} gives a smooth minimal surface $Y$
of general type with big cotangent bundle and
\[
 K_Y^2=117,\qquad c_2(Y)=207,\qquad s_2(Y)=-90.
\]
Let $Z$ be the canonical model of $Y$. Then $Z$ has eighteen
$A_8$ singularities and three $A_2$ singularities. Its natural
orbifold $\mathcal Z_{\rm can}$ satisfies
\[
 c_2(\mathcal Z_{\rm can})=39,\qquad
 \frac{c_1^2(\mathcal Z_{\rm can})}{c_2(\mathcal Z_{\rm can})}=3.
\]
Moreover, $Z$ is a compact complex two-ball quotient:
$
 Z\simeq\mathbb B^2/\Gamma,
$
where $\Gamma<\mathrm{PU}(2,1)$ is a discrete cocompact subgroup
acting freely in codimension one.
\end{theorem}

To the best of our knowledge, this is the first example of a
singular compact complex two-ball quotient whose smooth minimal
resolution has negative second Segre number and big cotangent
bundle. Earlier singular ball quotients include quotients of fake
projective planes~\cite{Keum08}.

Let us briefly explain the construction. A ramified base change of
a semi-stable fibration introduces chains of $(-2)$-curves above
the nodes of the fibers. Contracting specified chains gives a Du
Val quotient, and the QS criterion provides a numerical test for
the bigness of the cotangent bundle of its smooth resolution.
For Theorem~\ref{thm:sharp-three}, the initial genus-two fibration
already contains three $A_2$ chains. After the degree-nine base
change, these chains combine with the new exceptional curves.
Contracting all the resulting canonical configurations gives the
orbifold slope three. 

The proof of Theorem~\ref{thm:full-rational-interval} combines
three explicit semi-stable fibrations of genera two and three.
One construction produces orbifold Chern slopes arbitrarily close to
one. Genus-three interpolation and genus-two mixed-ramification
covers then yield overlapping rational intervals. Together with
the canonical-model interval obtained from the sharp example, these
intervals cover all of $(1,3]\cap\mathbb Q$. The bigness of the
cotangent bundle is established using the quotient-singularity
criterion of Asega--De Oliveira--Weiss. Our use of ramified base
changes is motivated by the methods of Liu--Lu and
Liu--L\"u~\cite{LiuJunLu,LiuXinLu}. The quotients in the lower
part of the range need not be canonical models; thus our orbifold
Chern slopes differ from the smooth Chern slopes studied
in~\cite{RU}.

\vspace{2mm}
The paper is organized as follows. Section~\ref{pre} reviews
surface fibrations, orbifold Chern numbers, base-change formulas,
bigness criteria and rational interpolation. In
Section~\ref{sec:sharp-endpoint} we prove
Theorem~\ref{thm:sharp-three}, establish ball quotient
uniformization, and construct canonical-model slopes approaching
three from below (Corollary~\ref{cor:tu-strict-near-three}).
Section~\ref{sec:aow-density} proves
Theorem~\ref{thm:full-rational-interval} by joining the genus-three
interpolation intervals to the genus-two constructions and the sharp
endpoint interval. Finally, Section~\ref{sec:further} proves an
optimal canonical-model slope for a genus-three pencil and gives
further asymptotic results and examples.

\section{Preliminaries and bigness criteria}\label{pre}
\subsection{Surface fibrations}
We recall some facts about surface fibrations and fix the notations used throughout.

Let $f: X\to C$ be a relatively minimal fibration from a smooth projective surface $X$ to a smooth projective curve $C$. Thus $f$ is a surjective morphism with connected fibers and no $(-1)$-curves in its fibers. Its general fiber is a smooth projective curve of genus $g\geq1$.
For a reduced connected nodal fiber $F$, let $\widetilde D$ be the normalization of an irreducible component $D$, and let $n_D$ count the branches of the nodes on $\widetilde D$ (a self-node contributes two). The fiber is semi-stable if $2g(\widetilde D)-2+n_D\geq0$ for every $D$, and stable if all these inequalities are strict. 
We say that $f$ is semi-stable if all its singular fibers are semi-stable.

Let $\omega_X$ (resp. $\omega_C$) be the canonical sheaf of $X$ (resp. $C$). Denote by $\omega_{X/C}=\omega_X\otimes f^*\omega_C^{-1}$
the relative canonical sheaf. Let $K_X$ (resp. $K_f$) be the divisor corresponding to $\omega_X$ (resp. $\omega_{X/C}$).
Let $b=g(C), q=h^1(X,\mathcal O_X)$, $q_f=q-b$, $p_g=h^0(X,\omega_X)$, $\chi(\mathcal O_X)=1-q+p_g$ and $\chi_{\text{top}}$
be the topological Euler characteristic of $X$.
The relative invariants of $f$ are defined by
\begin{equation*}
\begin{cases}
K_f^2   &= \omega_X^2-8(g-1)(b-1)=c_1^2(X)-8(g-1)(b-1),                 \\
 e_f     &= \chi_{\rm top}(X)-4(g-1)(b-1)=c_2(X)-4(g-1)(b-1),    \\
\chi_f &= \chi(\mathcal{O}_X)-(g-1)(b-1).
\end{cases}
\end{equation*}

These invariants are nonnegative and satisfy Noether's formula
\[12 \chi_f =K_f^2 +e_f.\]
For a semi-stable fibration, $e_f$ equals the total number of nodes in the singular fibers. 
For a relatively minimal non-isotrivial fibration $f$ of genus $g\geq2$, one has $\chi_f>0$, and its slope is defined by
\[\lambda_f=\frac{K_f^2}{\chi_f}.\]
We have $\lambda_f\leq12$. If equality holds, then $f$ is smooth; in this case, $f$ is called a Kodaira fibration.

\subsection{Du Val quotients and orbifold Chern numbers}\label{subsec:models}
We use the definition of an orbifold given in \cite[Section 2]{RR14}.
 Let $X$ be a normal algebraic variety and $\Delta =\sum d_jD_j$ be a $\mathbb{Q}$-divisor with $0 \leq d_j \leq 1$ on $X$ such that the log canonical divisor 
 $K_X+\Delta$ is $\mathbb{Q}$-Cartier. We call $(X, \Delta)$ a log pair.
 
 Now let $(X,\, \Delta)$ be a special log pair such that $\Delta=\sum_{j}(1-\frac{1}{m_j})D_j$, where the $D_j$ are prime divisors on $X$ and $m_j\in \mathbb{N}$. Then we define an orbifold as follows.
\begin{definition}
An orbifold $\mathcal{X}$ is a log pair $(X,\,\Delta=\sum_{j}(1-\frac{1}{m_j})D_j)$ such that $X$ is covered by orbifold charts $\varphi$ compatible with $\Delta$, i.e., $\varphi: U \to \varphi(U) \subset X$ is a Galois covering such that:
\begin{itemize}
    \item [(1)] $U$ is a domain in $\mathbb{C}^n$ and $\varphi(U)$ is open in $X$,
    \item [(2)] the divisorial part of the branch locus of $\varphi$ is $\operatorname{Supp}(\Delta)\cap\varphi(U)$; ramification over the singular locus of $X$ in codimension at least two is also allowed,
    \item [(3)] for any $x\in U'':=U\setminus \varphi^{-1}(X_{sing}\cup \Delta_{sing})$ such that $\varphi(x) \in D_j$, the ramification
index of $\varphi$ at $x$ satisfies $ord_{\varphi}(x) = m_j$.
\end{itemize}
\end{definition}

If $\Hat{X}$ has only ADE singularities (hence quotient singularities), the pair $(\Hat{X},0)$ carries its natural orbifold structure $\mathcal{X}$, 
denoted by $\phi: \mathcal{X} \to (\Hat{X}, \Delta)$.

Let $a_n$ (respectively, $d_n,\,e_n$) be the number of $A_n$
(respectively, $D_n,\,E_n$) singularities on $\Hat{X}$, and let $\Tilde{X}\to \Hat{X}$ be its minimal resolution. The orbifold Chern numbers are given by the following formula; see \cite[Proposition~13]{RR}.
\begin{pro}\label{obd}
    The Chern numbers of the orbifold $\mathcal{X}$ are $c_1^2(\mathcal{X})=c_1^2(\Hat{X})=c_1^2(\Tilde{X})$ and
    \begin{equation}\label{Ochern}
    \begin{aligned}
        c_2(\mathcal{X})={}&c_2(\Tilde{X})
        -\sum_{n\geq1}\left(n+1-\frac{1}{n+1}\right)a_n
        -\sum_{n\geq4}\left(n+1-\frac{1}{4(n-2)}\right)d_n\\
        &-\left(7-\frac{1}{24}\right)e_6
         -\left(8-\frac{1}{48}\right)e_7
         -\left(9-\frac{1}{120}\right)e_8.
    \end{aligned}
    \end{equation}
\end{pro}

We distinguish two orbifolds associated with a semi-stable pullback. An \emph{auxiliary} ADE orbifold $\mathcal X_{\rm aux}$ is defined relative to a designated initial semi-stable fibration: we contract the explicitly specified exceptional chains introduced by the chosen base changes, including their pullbacks and composite chains under later covers, while retaining the rational $(-2)$-components of the designated initial fibration.
  The \emph{canonical-model}
orbifold $\mathcal X_{\rm can}$ is obtained by contracting \emph{all}
$K$-trivial curves.  They need not have the same orbifold Euler number.
For the Roulleau--Rousseau criterion (Theorem~\ref{RR1}) we use the
canonical model.  For the QS-bigness criterion below, the auxiliary
quotient need not be canonical.

\subsection{Semi-stable base changes}
Let $f:S\to C$ be a surface fibration of genus $g\geq1$, and let $\pi:\Tilde C\to C$ be a finite morphism from a smooth projective curve. Form the fiber product $S\times_C\Tilde C$, with projection $\Pi':S\times_C\Tilde C\to S$. Let $\rho_1:S_1\to S\times_C\Tilde C$ be its normalization and $\rho_2:S_2\to S_1$ its minimal resolution. Finally, let $\Tilde\rho:S_2\to\Tilde S$ be the contraction to the relatively minimal model. These maps fit into the following commutative diagram:

\begin{center}
         \begin{tikzcd}
&\Tilde{S}\arrow[d,"\Tilde{f}"]\arrow[r, leftarrow,"\Tilde{\rho} "] & S_2 \arrow[d,"f_2 "] \arrow[r,"\rho_2 "] & S_1  \arrow[d,"f_1 "] \arrow[r,"\rho_1 "] & S\times_{C}\Tilde{C}  \arrow[d] \arrow[r,"\Pi' "] & S \arrow[d,"f"] \\
& \Tilde{C} \arrow[r,"id "]& \Tilde{C}\arrow[r,"id "]&\Tilde{C}\arrow[r,"id "]&\Tilde{C}\arrow[r,"\pi "]& C
\end{tikzcd}\,\,\,.
     \end{center}       

We call $\pi: \Tilde{C}\to C$ a base change and $\Tilde{f}$ the induced fibration. Its local behavior is determined in a neighborhood of each branch point of $C$ (see \cite{Xi90}).

Let $F_1,\ldots,F_s$ be the singular fibers of $f$, and let $\mathcal R=\{p_i=f(F_i):i=1,\ldots,s\}$ be the set of singular values.
The following local observation records the exceptional chains produced by the base change.
\begin{lemma}\label{bc}
Assume that $f$ is semi-stable, and let $\Tilde f:\Tilde X\to\Tilde C$ be the fibration induced by a degree-$d$ base change, with $d\geq2$, totally ramified over $p_1=f(F_1),\ldots,p_k=f(F_k)$, where $k\leq s$, and over $\ell$ additional points of $C\setminus\mathcal R$, and unramified elsewhere. If $N_0$ is the total number of nodes in $F_1,\ldots,F_k$, then the base change introduces $N_0(d-1)$ new $(-2)$-curves in the singular fibers of $\Tilde f$.
\end{lemma}
\begin{proof}
    It suffices to check the assertion at each node of the selected fibers.     
   At a node $q$ of $F_i$, choose local coordinates $(x,y)$ in which $f$ is given by $t=xy$. The totally ramified base change $t=u^d$ gives $xy=u^d$, an $A_{d-1}$ singularity. Its minimal resolution inserts a chain of $d-1$ new $(-2)$-curves in the fiber, as claimed.
\end{proof}

\begin{remark}
In Lemma~\ref{bc}, if the smooth pullback $\Tilde X$ is of general type, contracting only the new exceptional $(-2)$-chains gives a projective normal auxiliary surface $\Hat X$ of general type with precisely $N_0$ canonical quotient singularities of type $A_{d-1}$. It is not required to be the canonical model; the QS criterion of \cite{AOW0} applies to its minimal resolution. A sufficient condition for $\Tilde X$ to be of general type in the non-isotrivial constructions below is that $g\geq2$ and $g(\Tilde C)\geq2$.
\end{remark}

For a genus-$g$ semi-stable fibration over $\mathbb P^1$ with relative
invariants $K_f^2,e_f$, choose $k$ singular fibers containing $N_0$
nodes and $\ell$ smooth fibers, and put $m=k+\ell$.  Suppose that a
connected degree-$d$ cover is totally ramified at exactly these $m$
points and nowhere else.  Set
\[
V=d(m-2)-m,\qquad
A=dK_f^2+4(g-1)V,\qquad E=de_f+2(g-1)V.
\]
Then, for the smooth semi-stable pullback $Y$,
\begin{equation}\label{eq:AE-general}
c_1^2(Y)=A,\quad c_2(Y)=E,\quad s_2(Y)=A-E,
\qquad 2g(C_d)-2=V.
\end{equation}
The auxiliary orbifold has
\begin{equation}\label{eq:aux-general}
B_{\rm aux}:=c_2(\mathcal X_{\rm aux})
=E-N_0\left(d-\frac1d\right).
\end{equation}
Such a cover exists whenever $d\geq2$, $m\geq3$ and
$(d-1)m$ is even; the branching cycles can be chosen as inverse
pairs, together with a triple of $d$-cycles when both $d,m$ are odd.

\begin{lemma}\label{S2}\leavevmode
Let $f:X\to\mathbb P^1$ be a relatively minimal semi-stable fibration of genus $g\geq2$. Take a connected degree-$d$ base change, with $d\geq2$, totally ramified over $k$ singular values whose fibers contain $N_0$ nodes in total and over $\ell$ smooth-fiber points, and unramified elsewhere. Let $\mathcal X$ denote the auxiliary orbifold obtained by contracting the new chains. Then
    \begin{itemize}
        \item[(I)] 
        \begin{equation}\label{Segre1}
        s_2(\Tilde{X})=d(K_f^2-e_f)+2(g-1)(d(k+\ell-2)-k-\ell),
    \end{equation}

    \item[(II)] \begin{equation}\label{Segre2}
       s_2(\mathcal{X})=d\cdot (K_f^2-e_f+N_0)-\frac{N_0}{d}+2(g-1)(d(k+\ell-2)-k-\ell),
    \end{equation}

 \item[(III)]
     \begin{equation}\label{Segre3}
     s_2(\Tilde{X})+s_2(\mathcal{X})=d(2K_f^2-2e_f+N_0)-\frac{N_0}{d}+4(g-1)(d(k+\ell-2)-k-\ell).
 \end{equation}
   \end{itemize}
    
\end{lemma}
\begin{proof}
The semi-stable base-change formulas and Riemann--Hurwitz give
\begin{equation*}
 K_{\widetilde f}^{2}=dK_f^2,\qquad e_{\widetilde f}=de_f,\qquad
 2(\widetilde b-1)=d(k+\ell-2)-k-\ell.
\end{equation*}
Since
$c_1^2(\widetilde X)=K_{\widetilde f}^{2}+8(g-1)(\widetilde b-1)$
and
$c_2(\widetilde X)=e_{\widetilde f}+4(g-1)(\widetilde b-1)$,
subtraction gives~\eqref{Segre1}.  The $N_0$ contracted
$A_{d-1}$ chains have total orbifold Euler correction
$N_0(d-1/d)$ by~\eqref{Ochern}; hence
\[
 c_1^2(\mathcal X)=c_1^2(\widetilde X),\qquad
 c_2(\mathcal X)=c_2(\widetilde X)-N_0\left(d-\frac1d\right).
\]
Thus $s_2(\mathcal X)=s_2(\widetilde X)+N_0(d-1/d)$,
which proves~\eqref{Segre2}; adding~\eqref{Segre1}
gives~\eqref{Segre3}.
\end{proof}

\subsection{Bigness criteria}
\label{subsec:qs}
We first recall the criterion of Roulleau and Rousseau
\cite[Theorem~1]{RR14} for a \emph{canonical} surface.
\begin{theorem}\label{RR1}
Let $\widehat X$ be a projective surface with ample canonical divisor
and at worst canonical (Du Val) singularities, let
$Y\to\widehat X$ be its minimal resolution, and let $\mathcal X$
be the orbifold associated with $\widehat X$.  If
\begin{equation*}
 s_2(Y)+s_2(\mathcal X)>0,
\end{equation*}
then $\Omega_Y^1$ is big.  In particular, $Y$ has only finitely
many rational or elliptic curves.
\end{theorem}
The original nodal-hypersurface argument of
Bogomolov--De Oliveira \cite{BD} motivates this approach, but its
numerical claim requires the correction discussed in
\cite[Section~4.1]{RR14}.  In applying the preceding theorem we
therefore use its stated canonical-model hypotheses.

Asega--De Oliveira--Weiss instead give the QS-bigness criterion in terms of
\emph{individual quotient singularities} which need not lie on the
canonical model.  This distinction is essential for our auxiliary
construction.

In the auxiliary constructions below, $Y$ denotes the \emph{smooth minimal resolution}
of the \emph{auxiliary quotient surface} $Z_{\mathrm{aux}}$,
obtained by contracting the specified new exceptional chains introduced by
the ramified base changes.  In particular, the old $(-2)$-bridges of the
quartic pencil are \emph{not} contracted.  This convention applies to auxiliary constructions. In constructions explicitly using the canonical model, all $K$-trivial curves, including the old bridges, are contracted.
The following criterion does not require $Z_{\mathrm{aux}}$ to be the
canonical model.

Write $h_n=h^1_\Omega(A_n)$ for the local first-cohomology
asymptotics determined by Asega--De Oliveira--Weiss
\cite{AOW0,AOW}.  Their QS-bigness criterion
\cite[Theorem~1.1]{AOW0} says that if a smooth surface $Y$ of general
type is the minimal resolution of a projective normal surface with
quotient singularities $z_j$, then
\begin{equation}\label{eq:qs-criterion}
 s_2(Y)+6\sum_j h^1_\Omega(z_j)>0
 \quad\Longrightarrow\quad \Omega^1_Y\text{ is big}.
\end{equation}
For $n\geq1$, Asega--De Oliveira--Weiss also prove
\cite[Corollary~2.11]{AOW0}
\[
 \frac n6-\frac{4\pi^2-39}{18}<h_n<\frac n6.
\]
Because $(4\pi^2-39)/3<1/5$, we have the useful rational bound
\begin{equation}\label{eq:qs-local-bound}
 6h_{d-1}>d-\frac65\qquad(d\geq2).
\end{equation}

\begin{proposition}
\label{prop:aow-certificate}
Let $f$ be a genus-$g$ semi-stable fibration over $\mathbb P^1$, with $g\geq2$, and let $Y$ be the smooth semi-stable pullback under a connected degree-$d$ cover, with $d\geq2$, totally ramified over the points corresponding to $k$ selected singular fibers and $\ell$ smooth fibers, and unramified elsewhere. Let $N_0>0$ be the total number of nodes in the selected singular fibers, and put
\[
 V=d(k+\ell-2)-(k+\ell),\qquad
 S=d(K_f^2-e_f)+2(g-1)V.
\]
Assume $Y$ is a minimal surface of general type and only the new
$N_0(d-1)$ exceptional $(-2)$-curves are contracted.
Then $Z_{\mathrm{aux}}$ has precisely $N_0$ singularities of type
$A_{d-1}$ and $s_2(Y)=S$.  In particular,
\begin{equation*}
  S+N_0\left(d-\frac65\right)\geq0
  \quad\Longrightarrow\quad\Omega^1_Y\text{ is big}.
\end{equation*}
This also holds when $Y$ has old $(-2)$-curves which are not contracted
by $Y\to Z_{\mathrm{aux}}$.
\end{proposition}
\begin{proof}
The local equation at each selected node is $xy=t$.  Under total
ramification $t=u^d$, it becomes $xy=u^d$; the new exceptional chain
has type $A_{d-1}$.  There is one such chain per selected original
node.  All the old curves, if any, remain on $Y$ and introduce no
additional singularities on $Z_{\mathrm{aux}}$.
Since $Y$ is the minimal resolution of $Z_{\mathrm{aux}}$,
\eqref{eq:qs-criterion} and \eqref{eq:qs-local-bound} imply the claim.
\end{proof}

\begin{remark}
\label{rem:aow-versus-rr}
Put $B=c_2(\mathcal X_{\mathrm{aux}})$ and $A=c_1^2(Y)$,
$E=c_2(Y)$.  The auxiliary Segre-sum condition is
$2A-E-B>0$, equivalently
$2S+N_0(d-1/d)>0$.  For $d\geq2$,
\[
 d-\frac65\ \geq\ \frac12\left(d-\frac1d\right),
\]
so this strict Segre-sum inequality implies the QS condition.
The converse fails in general; thus the QS criterion provides
additional examples even when the auxiliary quotient is not canonical.
The auxiliary orbifold Chern slope here is a numerical invariant
of $Z_{\mathrm{aux}}$, not the fibration slope $\lambda_f$.
\end{remark}

\begin{proposition}\label{prop:orbifold-slope-range}
Let $Z$ be a projective normal surface with only Du Val singularities,
and suppose that its minimal resolution $Y$ is a minimal surface of
general type. Let $\mathcal Z$ be its natural orbifold. Then
\[
 c_2(\mathcal Z)>0,\qquad
 0<\frac{c_1^2(\mathcal Z)}{c_2(\mathcal Z)}\leq3.
\]
Moreover, $s_2(\mathcal Z)>0$ if and only if this slope is greater
than one. If all the singularities are of type $A_n$ and the QS
inequality~\eqref{eq:qs-criterion} holds, then $s_2(\mathcal Z)>0$.
The same conclusion holds when $s_2(Y)<0$ and
$s_2(Y)+s_2(\mathcal Z)>0$.
\end{proposition}
\begin{proof}
Write $\pi:Y\to Z$ for the resolution. It is crepant, so
$K_Y=\pi^*K_Z$. Thus $K_Z$ is nef
and big. The orbifold Bogomolov--Miyaoka--Yau inequality
\cite{Miy84,GKPT19} gives
\[
 0<K_Y^2=c_1^2(\mathcal Z)\leq3c_2(\mathcal Z).
\]
This proves the first assertion. The equivalence follows from
the definition of $s_2(\mathcal Z)$.

For a singularity of type $A_n$, let
$\delta_n=n+1-1/(n+1)$ be its orbifold Euler correction.
By~\cite[Corollary~2.11]{AOW0}, we have
$6h_n<n<\delta_n$. Hence, if the singularities have types
$A_{n_1},\ldots,A_{n_q}$ and the QS inequality holds, then
\[
 s_2(\mathcal Z)
 =s_2(Y)+\sum_{j=1}^q\delta_{n_j}
 \geq s_2(Y)+6\sum_{j=1}^q h_{n_j}>0.
\]
Finally, the Segre-sum hypothesis with $s_2(Y)<0$ gives
$s_2(\mathcal Z)>-s_2(Y)>0$.
\end{proof}

\subsection{Interpolation over smooth fibers}\label{subsec:qs-interpolation}
We shall use the following interpolation for both auxiliary
quotients and canonical models. The singularities may have
different ADE types.
\begin{proposition}
\label{prop:aow-interp}
Let $Y\to C$ be a relatively minimal semi-stable genus-$g$ fibration,
$g\geq2$, with $g(C)\geq2$. Let $Y\to Z$ contract a specified union
of vertical ADE configurations, and let $\mathcal X$ be the natural
orbifold of $Z$. Put
\[
 A=c_1^2(Y),\quad E=c_2(Y),\quad B=c_2(\mathcal X),\quad
 H=\sum_{z\in\operatorname{Sing}Z}h^1_\Omega(z).
\]
Assume $0<A<E$, $B>0$, and $A-E+6H>0$.
For every rational $t\in[0,(E-A)/2)$, a further connected base change
branched only above smooth fibers yields a smooth minimal surface
$Y_t$ and a Du Val quotient with orbifold $\mathcal X_t$ such that
\[
 s_2(Y_t)<0,\qquad\Omega^1_{Y_t}\text{ is big},\qquad
 \frac{c_1^2(\mathcal X_t)}{c_2(\mathcal X_t)}
 =\frac{A+4t}{B+2t}.
\]
The resulting set of rational slopes is exactly the rational points
in the interval with endpoints
\[
 r_0=\frac AB\quad\text{(included)},\qquad
 r_1=\frac{2E-A}{B+E-A}\quad\text{(excluded)},
\]
provided $A\ne2B$; if $A=2B$, the slope is constantly $2$.
If $Z$ is the canonical model of $Y$, the resulting quotients are
also canonical models.
\end{proposition}
\begin{proof}
Take a degree-$n$ cover $C'\to C$ with total ramification degree
$R$, supported at points where the fibration is smooth. Set
$t=(g-1)R/n$. Riemann--Hurwitz and semi-stable base change give
\[
 \frac{c_1^2(Y')}{n}=A+4t,\qquad
 \frac{c_2(Y')}{n}=E+2t.
\]
The cover is unramified above every singular fiber, so each
contracted ADE configuration lifts to $n$ disjoint copies of the
same type. Contracting these copies yields a projective Du Val
quotient with
\[
 \frac{c_2(\mathcal X')}{n}=B+2t,\qquad H'=nH.
\]
Thus
\[
 \frac{s_2(Y')}{n}=A-E+2t<0,\qquad
 \frac{s_2(Y')+6H'}n=A-E+2t+6H>0.
\]
The QS criterion~\eqref{eq:qs-criterion} proves bigness.
The pullback is relatively minimal, and its base has genus at least
two; hence it is a minimal surface of general type.

Every rational $t>0$ can be written as $t=2(g-1)h/n$ with
$n\geq2$, $h\geq1$. Choose $2h$ distinct smooth-fiber points on
$C$. Send one handle generator to an $n$-cycle, its partner to the
identity, and the small loops about these points to identical
transpositions in pairs. These permutations satisfy the
surface-group relation and generate a transitive subgroup, so the
Riemann existence theorem gives the required connected cover.
For $t=0$, take the identity. Finally,
\[
 \frac{d}{dt}\frac{A+4t}{B+2t}
 =\frac{4B-2A}{(B+2t)^2},
\]
which gives the interval and its endpoint conventions.
If $Z$ is canonical, all $K_Y$-trivial curves are vertical. Their
copies exhaust the $K_{Y'}$-trivial curves: branching over smooth
fibers creates no exceptional curves, and the cover is unramified
over all singular fibers. Their contraction is therefore the
canonical model of $Y'$.
\end{proof}

\begin{corollary}\label{prop:interp}
Let $Y\to C$ be a relatively minimal semi-stable fibration of genus
$g\geq2$ over a curve with $g(C)\geq2$, and let $Z$ be its canonical
model with natural orbifold $\mathcal X_{\rm can}$. Put
\[
 A=c_1^2(Y),\qquad E=c_2(Y),\qquad B=c_2(\mathcal X_{\rm can}),
\]
and assume $0<A<E$, $B>0$ and $2A-E-B>0$.
For every rational $t\in[0,(E-A)/2)$, a further connected base change
branched only above smooth fibers gives a smooth minimal surface $Y_t$
whose canonical-model orbifold $\mathcal X_t$ satisfies
\begin{equation}\label{eq:interp}
 \frac{c_1^2(\mathcal X_t)}{c_2(\mathcal X_t)}
   =\frac{A+4t}{B+2t},\qquad
 \frac{s_2(Y_t)}n=A-E+2t,\qquad
 \frac{s_2(Y_t)+s_2(\mathcal X_t)}n=2A-E-B+4t,
\end{equation}
where $n$ is the degree of the additional cover. In particular,
$s_2(Y_t)<0$ and $\Omega^1_{Y_t}$ is big. These slopes comprise the
rational values between $A/B$ (included) and
$(2E-A)/(B+E-A)$ (excluded), unless $A=2B$, when the slope is
constantly $2$.
\end{corollary}
\begin{proof}
Use the connected covers constructed in the proof of
Proposition~\ref{prop:aow-interp}, with $t=(g-1)R/n$ and total
ramification $R$ supported on smooth fibers. They are unramified
above the singular fibers, so the canonical singularities lift in
$n$ copies and the resulting quotient remains the canonical model.
The base-change formulas give~\eqref{eq:interp}. The assumed Segre
sum is positive and increases by $4t$, while $s_2(Y_t)<0$ in the
stated range. Theorem~\ref{RR1} proves bigness. Finally, the
fractional-linear slope is strictly monotone unless $A=2B$.
\end{proof}

\section{The sharp endpoint and ball quotient uniformization}\label{sec:sharp-endpoint}\label{Sec 3}\label{sec:geography}
\subsection{The initial genus-two fibration}\label{subsec:tu-five}
\label{sub1}
We start with the genus-two fibration constructed in
\cite[Theorem~1.1(2) and Section~3]{Tu07}.
Let $R_1$ be the curve defined by $x=y^3$ on $P_1=\mathbb{P}^1 \times \mathbb{P}^1$. We construct a double covering $\pi_1: P_2 \to P_1 $, branched over the curves $x = 0$ and $x = 1$. Let $R_2 = \pi_1^{*}(R_1)$ and let $ \pi_2: P_3 \to P_2$ be a
triple cover totally ramified over the pullback of $x = \infty$ on $ P_2$. Then $P_3=\mathbb{P}^1 \times \mathbb{P}^1$.
Let $pr_2$ denote the projection from $P_3$ to the second factor, let $F_{\infty}$ be the fiber of $pr_2$ 
over $y = \infty$, and let $R_3 = \pi_2^*(R_2)$. Then $R = R_3 + F_{\infty} $ is an even divisor of type $(6, 4)$.
Consider the double cover of $P_3$ branched over $R$. Resolving its singularities and contracting the vertical $(-1)$-curve arising from $F_\infty$ yields a semi-stable fibration $f: X \to \mathbb{P}^1 $  of genus $2$ with $5$ singular fibers (the fiber over $y = \infty$ is the union of two smooth elliptic curves).
The invariants can be computed using the standard formulas for double covers (cf. \cite{BPV}) as follows; see \cite[Section~3]{Tu07} for more details:
\[p_g(X) = q(X) = 0, K^2_f = 5, e_f=19, \chi_f=2.\]

\begin{lemma}\label{lem:tu-canonical-chains}
For the five-fiber genus-two fibration of Tu in
Subsection~\ref{subsec:tu-five}, label the singular values by
$p_0=0$, $p_1=1$, $p_2=\omega$, $p_3=\omega^2$, and $p_4=\infty$,
where $\omega$ is a primitive cube root of unity. The fiber over $p_0$
contains precisely three old $A_2$ chains. After contracting them, its
stable model consists of two smooth rational curves meeting at three
nodes. The fibers over $p_1,p_2,p_3$ are already stable and each consists
of two smooth rational curves meeting at three nodes. The fiber over
$p_4$ consists of two smooth elliptic curves meeting at one node.
In particular, the semi-stable node counts are $(9,3,3,3,1)$.

Let $C'\to\mathbb P^1$ be a connected finite cover of degree $n$ with
$g(C')\geq2$, and let $e_{ij}$ be the local degrees above $p_i$, so
$\sum_j e_{ij}=n$. On the smooth semi-stable pullback $Y$, the
canonical-model chains are three chains of type $A_{3e_{0j}-1}$ over
each point above $p_0$, three chains of type $A_{e_{ij}-1}$ over each
point above $p_i$, $i=1,2,3$, and one chain of type $A_{e_{4j}-1}$
over each point above $p_4$. Here a chain of type $A_0$ is empty.
These are all the $K_Y$-trivial curves.
\end{lemma}

\begin{proof}
Choose coordinates on the two successive ruled-surface covers in Tu's
construction so that
\[
 x=\frac{u^2}{u^2-1},\qquad
 u=\frac{z^3+1}{z^3-1},\qquad
 x=\frac{(z^3+1)^2}{4z^3}.
\]
The general fiber of the double cover then has affine equation
\[
 w^2=z^6+(2-4y^3)z^3+1.
\]
At $y=0$, the two components $w=\pm(z^3+1)$ meet at the three roots
of $z^3+1$. The total space has local equation $ab=y^3$ up to
analytic units at each of these points. Resolving these three $A_2$
singularities inserts three chains of two $(-2)$-curves and gives
nine nodes. At each value $y^3=1$, the fiber is
$w^2=(z^3-1)^2$; its two rational components meet at three nodes,
and the total space is smooth there. Tu's construction identifies the
remaining fiber, over infinity, as the union of two smooth elliptic
curves; see~\cite[Theorem~1.1(2) and Section~3]{Tu07}. Since it is
semi-stable of genus two, the two components meet at exactly one node.

The stable smoothing equation at a node over $p_0$ is $ab=y^3$.
A base change of local degree $e$ gives $ab=v^{3e}$, hence a chain
of type $A_{3e-1}$. At the other stable nodes the equation is
$ab=t$, and local degree $e$ gives a chain of type $A_{e-1}$.

For each local degree $e$ above $p_0$, the canonical correction is
$3(3e-1/(3e))=9e-1/e$, whereas contracting only the new chains
would contribute $9(e-1/e)$. Consequently,
\[
 B_{\rm can}=B_{\rm aux}-8\sum_j\frac1{e_{0j}}.
\]

Here $B_{\rm can}$ and $B_{\rm aux}$ are the corresponding orbifold
second Chern numbers.

The smooth pullback is relatively minimal. Since the base has genus
at least two, every rational curve on $Y$ is vertical, so $Y$ is
minimal. It is of general type by the standard positivity of the
relative canonical divisor for a genus-two fibration over such a base.
If $K_Y\cdot D=0$ for an irreducible curve $D$, the Hodge index
theorem and adjunction show that $D$ is a smooth rational
$(-2)$-curve. Each surviving rational component of a stable fiber
meets three branches and has canonical degree one. Each elliptic
component over infinity also has canonical degree one. The listed
chains therefore exhaust the $K_Y$-trivial curves, and their
contraction is the canonical model.
\end{proof}

\subsection{Proof of Theorem~\ref{thm:sharp-three}}
\begin{proof}
Use the five-fiber genus-two fibration of
Lemma~\ref{lem:tu-canonical-chains}, with relative invariants
$K_f^2=5$ and $e_f=19$. Let $C$ be the smooth projective
normalization of
\[
 u^9=y^3(y-1)(y-\omega)(y-\omega^2)^7.
\]
The projection $C\to\mathbb P^1_y$ is a connected cyclic cover of
degree nine: the valuation at $y=1$ is one, so its defining rational
function is not a proper power corresponding to a divisor of nine.
The valuations at the five singular values are
$3,1,1,7,-12$. Consequently the local ramification indices are
\[
 (3,9,9,9,3),
\]
with three points above each value of index three and one point
above each value of index nine. There is no other branching.
Riemann--Hurwitz gives
\[
 R=2\cdot3(3-1)+3(9-1)=36,\qquad
 2g(C)-2=-18+36=18.
\]
Thus $g(C)=10$. For the smooth semi-stable pullback $Y$,
\[
 \begin{aligned}
 c_1^2(Y)&=9\cdot5+8(10-1)=117,\\
 c_2(Y)&=9\cdot19+4(10-1)=207.
 \end{aligned}
\]
The surface is minimal of general type by
Lemma~\ref{lem:tu-canonical-chains}. That lemma gives nine $A_8$
chains above $p_0$, nine further $A_8$ chains above $p_1,p_2,p_3$,
and three $A_2$ chains above $p_4$. Contracting them gives the
canonical model and hence
\[
 c_2(\mathcal Z)
 =207-18\left(9-\frac19\right)
      -3\left(3-\frac13\right)=39.
\]
The resolution is crepant, so $c_1^2(\mathcal Z)=117$.

It remains to prove bigness of the cotangent bundle. The rational
local estimate~\eqref{eq:qs-local-bound} gives
\[
 \begin{aligned}
 s_2(Y)+6\sum_{z\in\operatorname{Sing}Z}h^1_\Omega(z)
 &>-90+18\left(9-\frac65\right)
          +3\left(3-\frac65\right)\\
 &=\frac{279}{5}>0.
 \end{aligned}
\]
The QS criterion~\eqref{eq:qs-criterion} therefore implies that
$\Omega_Y^1$ is big.

Finally, $Z$ has Du Val singularities and hence is klt. Since $Z$
is the canonical model of $Y$, its canonical divisor $K_Z$ is ample.
The orbifold Chern numbers of $Z$ agree with its $\mathbb Q$-Chern
numbers, so the above computation gives
\[
 \widehat c_1^2(\mathscr T_Z)=K_Z^2=117
 =3\cdot39=3\widehat c_2(\mathscr T_Z).
\]
By the characterization of singular ball quotients
\cite[Theorem~1.8]{GKPT20} or the orbifold uniformization \cite{CGG24}, we have
\[
 Z\simeq\mathbb B^2/\Gamma
\]
for a discrete cocompact subgroup
$\Gamma<\mathrm{PU}(2,1)$ acting freely in codimension one.
\end{proof}

\subsection{Local ramification and further covers}
\begin{proposition}\label{prop:tu-optimal-square}
For any connected finite base change of the five-fiber genus-two
fibration whose new base has genus at least two, put
$A=c_1^2(Y)$ and $B=c_2(\mathcal Z)$, where $\mathcal Z$ is the
canonical-model orbifold. Let $R_{\rm sm}$ denote the total
ramification over smooth-fiber values, and retain the local degrees
$e_{ij}$ of Lemma~\ref{lem:tu-canonical-chains}. Then
\[
 3B-A
 =2R_{\rm sm}
  +\sum_{i\in\{0,4\}}\sum_j
       \frac{(e_{ij}-3)^2}{3e_{ij}}
  +\sum_{i=1}^{3}\sum_j
       \frac{(e_{ij}-9)^2}{9e_{ij}}.
\]
In particular, when $B>0$ the slope is at most three. Equality
holds precisely when there is no branching over smooth fibers,
all local degrees above $p_0,p_4$ are three, and all local degrees
above $p_1,p_2,p_3$ are nine. The cover in
Theorem~\ref{thm:sharp-three} realizes these conditions.
\end{proposition}

\begin{proof}
Let the covering degree be $n$, put $q=\sum_{i=0}^4\#\{j\}$,
and let $R=5n-q+R_{\rm sm}$. Semi-stable base change gives
\[
 A=4R-3n,\qquad E:=c_2(Y)=2R+15n.
\]
Using the complete canonical-chain list gives
\[
 B=2R-4n
   +\sum_j\frac1{e_{0j}}+\sum_j\frac1{e_{4j}}
   +3\sum_{i=1}^3\sum_j\frac1{e_{ij}}.
\]
Thus $3B-A=2R-9n+3\sum_{i,j}c_i/e_{ij}$ with
$(c_0,\ldots,c_4)=(1,3,3,3,1)$. Expand the squares and use
$\sum_j e_{ij}=n$ and $\sum_i1/(3c_i)=1$ to obtain the identity.
The equality assertion follows because all terms are nonnegative.
\end{proof}

\begin{corollary}\label{cor:tu-strict-near-three}
For every integer $m\geq1$, there is a smooth minimal surface $Y_m$
of general type with $s_2(Y_m)<0$ and big cotangent bundle whose
canonical-model orbifold $\mathcal Z_m$ has slope
\[
 \frac{c_1^2(\mathcal Z_m)}{c_2(\mathcal Z_m)}
 =3-\frac{2}{39m+2}.
\]
In particular, these slopes tend to three from below.
\end{corollary}
\begin{proof}
Start with $Y\to C$ from Theorem~\ref{thm:sharp-three}, where
$g(C)=10$. Choose two distinct points of $C$ supporting smooth
fibers. There exists a connected degree-$2m$ cover $C_m\to C$
simply branched at exactly these two points. Indeed, send one
handle generator to a $2m$-cycle, its paired generator to the
identity, the two peripheral loops to the same transposition,
and the remaining handle generators to the identity. The
surface-group relation is satisfied, the generated action is
transitive, and the Riemann existence theorem supplies the cover.
Its total ramification is two.

There are no new exceptional chains above smooth fibers, and the
canonical singularities lift to $2m$ disjoint copies. The smooth
pullback therefore satisfies
\[
 \begin{aligned}
 c_1^2(Y_m)&=234m+8,& c_2(Y_m)&=414m+4,\\
 c_2(\mathcal Z_m)&=78m+4,& s_2(Y_m)&=-180m+4<0.
 \end{aligned}
\]
The QS lower bound scales by $2m$, while smooth-fiber ramification
adds four to the smooth Segre number. It remains strictly positive,
so $\Omega_{Y_m}^1$ is big. The slope formula follows. For example, $m=1000$ gives the slope
$58502/19501\approx2.99994872$.
\end{proof}

\begin{corollary}\label{cor:sharp-three-interval}
Every rational slope in $(99/43,3]$ is attained by a canonical-model
orbifold whose smooth minimal resolution has
negative second Segre number and big cotangent bundle.
\end{corollary}
\begin{proof}
Apply the QS interpolation of Proposition~\ref{prop:aow-interp}
to the surface of Theorem~\ref{thm:sharp-three}, with
$(A,E,B)=(117,207,39)$. Its strict QS inequality was verified
above. The admissible rational parameters are $0\leq t<45$, and
\[
 r(t)=\frac{117+4t}{39+2t}
\]
is strictly decreasing, with $r(0)=3$ and
$\lim_{t\to45^-}r(t)=99/43$. For every rational $r$ in the asserted
interval, the inverse parameter is
\[
 t=\frac{39(3-r)}{2(r-2)}\in[0,45)\cap\mathbb Q.
\]
The proposition applies to singularities of different ADE types and
preserves canonical models, so each such rational slope is attained.
\end{proof}

\section{Realization of all rational orbifold Chern slopes}\label{sec:aow-density}
In this section we prove Theorem~\ref{thm:full-rational-interval}.
The upper interval was obtained in
Corollary~\ref{cor:sharp-three-interval}. We construct the lower
intervals using auxiliary quotients and then join them to this
canonical-model interval.
\subsection{Slopes approaching one}
\label{subsec:near-one}
We use the genus-two fibration of Subsection~\ref{subsec:tu-five} with five
singular fibers, constructed by Tu~\cite{Tu07}.  Its relative
invariants are $K_f^2=5$ and $e_f=19$.  One singular fiber
$F_\infty$ is the union of two smooth elliptic curves and has exactly
one node: its arithmetic genus is $1+1+\delta-1=2$, so $\delta=1$.
Allowing partial rather than only total ramification over smooth
fibers gives the following family.

\begin{proposition}\label{prop:tu-near-one}
Every rational number in $(1,33/32)$ is the orbifold Chern slope of
an auxiliary quotient with Du Val singularities whose smooth minimal
resolution is a minimal surface $Y$ of general type with
$s_2(Y)<0$ and big cotangent bundle.
\end{proposition}
\begin{proof}
Fix an integer $h\geq2$, and put $n=4h$ and $a=h+6\leq n$.
Choose nine distinct points of $\mathbb P^1$ over which the fibration
is smooth.  We construct a connected degree-$n$ cover
$\pi_h\colon C_h\to\mathbb P^1$ totally ramified at the point
corresponding to $F_\infty$ and at seven of the chosen smooth-fiber
points, and with cycle type $(a,1^{n-a})$ at the remaining two points.
Indeed, let $\sigma\in S_n$ be an $n$-cycle and $\tau\in S_n$ an
$a$-cycle.  The local monodromies
\[
 \sigma,\sigma^{-1},\sigma,\sigma^{-1},
 \sigma,\sigma^{-1},\sigma,\sigma^{-1},\tau,\tau^{-1}
\]
have product one and generate a transitive subgroup, since $\sigma$
is an $n$-cycle.  The Riemann existence theorem supplies the desired
cover.  The total ramification and the genus of $C_h$ are
\begin{equation}\label{eq:tu-near-one-genus}
 R_h=8(n-1)+2(a-1)=34h+2,\qquad
 2g(C_h)-2=R_h-2n=26h+2.
\end{equation}

Let $Y_h$ be the smooth semi-stable pullback.  The relative invariants
scale by $n$; using~\eqref{eq:tu-near-one-genus} and the standard
base-change formulas gives
\begin{equation*}
 \begin{aligned}
 c_1^2(Y_h)&=nK_f^2+4(R_h-2n)=124h+8,\\
 c_2(Y_h)&=ne_f+2(R_h-2n)=128h+4,\\
 s_2(Y_h)&=4-4h<0.
 \end{aligned}
\end{equation*}
The base curve has genus $13h+2\geq2$. The pulled-back fibration
is relatively minimal, and every rational curve is vertical; hence
$Y_h$ is minimal. Since both the base and the general fiber have
genus at least two, $Y_h$ is of general type.

Over the node of $F_\infty$, total ramification yields one
$A_{n-1}$ chain, by the local equation $xy=u^n$.
Contracting just this chain gives an auxiliary normal surface with
one $A_{n-1}$ singularity and orbifold Euler number
\begin{equation*}
 B_h:=c_2(\mathcal X_{h,\mathrm{aux}})
 =c_2(Y_h)-\left(n-\frac1n\right)
 =124h+4+\frac1{4h}.
\end{equation*}
By the local estimate~\eqref{eq:qs-local-bound},
\[
 s_2(Y_h)+6h^1_\Omega(A_{4h-1})
 >(4-4h)+4h-\frac65=\frac{14}{5}>0.
\]
Thus the QS criterion proves that $\Omega_{Y_h}^1$ is big, and
\[
 r_h:=\frac{c_1^2(Y_h)}{B_h}
 =1+\frac{4-1/(4h)}{124h+4+1/(4h)}\longrightarrow1^+.
\]

Finally apply the rational interpolation of
Proposition~\ref{prop:aow-interp} to $Y_h\to C_h$.  Since
$\bigl(c_2(Y_h)-c_1^2(Y_h)\bigr)/2=2h-2$, every rational
$t\in[0,2h-2)$ is admissible, and the corresponding orbifold
slopes fill
\begin{equation*}
 \left[r_h,\frac{2c_2(Y_h)-c_1^2(Y_h)}
 {B_h+c_2(Y_h)-c_1^2(Y_h)}\right)\cap\mathbb Q
 =\left[r_h,\frac{528h^2}{512h^2+1}\right)\cap\mathbb Q.
\end{equation*}
The left endpoints tend to $1$ from above, whereas the right
endpoints tend to $33/32$ from below.  Every rational number in
$(1,33/32)$ belongs to one of these intervals, as claimed.
\end{proof}
\begin{remark}
The quotient in this construction is not asserted to be the canonical
model.  Indeed, the original five-fiber genus-two fibration has
$19$ nodes, one of them on the stable fiber $F_\infty$.
A stable genus-two fiber has at most three nodes, so among the other
four fibers at least six rational $(-2)$-components must be contracted
in the relative stable model.  Our covers are unramified over these
four fibers.  The resulting older exceptional curves therefore also
contribute to the canonical-model orbifold Euler correction, and
cannot be ignored in the near-one calculation.
\end{remark}

\subsection{The genus-three fibration and its interpolation intervals}\label{sub2}
\label{subsec:four-g3}
We use a non-hyperelliptic genus-three fibration over $\mathbb P^1$.

Consider a pencil defined by $F+t\cdot G=0$ on the projective plane $\mathbb{P}^2$, where 
\[F=X^4 + Y^4 + Z^4,\,G=X^2Y^2 + X^2Z^2 + Y^2Z^2,\,t\in \mathbb{P}^1.\]
Define the rational map $\mathbb P^2\dashrightarrow\mathbb P^1$ by $[X,Y,Z]\dashrightarrow[-F:G]$, so that the affine coordinate is $t=-F/G$.

After blowing up the base points of the pencil, we obtain a non-hyperelliptic genus-three fibration $f_0: X_0 \to \mathbb{P}^1$ whose singular fibers are $F_2,\,F_{-2},\, F_{-1},\, F_{\infty}$ over $2,-2,-1,\infty$ respectively. These four singular fibers have the following forms; see also \cite[Section~5.2]{MR25} for the corresponding plane curves:

\begin{itemize}
    \item [(1)] $F_{-2}=\sum\limits_{i=1}^4\Gamma_i$ is stable, where the $\Gamma_i$ are $(-3)$-curves with $\Gamma_i\cdot\Gamma_j=1$ for $i\ne j$.
\item [(2)] $F_{2}=2\Gamma_0+\sum\limits_{i=1}^8\Gamma_i$, where $\Gamma_0$ is a $(-4)$-curve, the $\Gamma_i$ for $1\leq i\leq8$ are $(-2)$-curves, $\Gamma_0\cdot\Gamma_j=1$ for $1\leq j\leq8$, and $\Gamma_i\cdot\Gamma_j=0$ for $1\leq i\ne j\leq8$.
\item [(3)] $F_{-1}=\Gamma_1+\Gamma_2$ is stable,  where $\Gamma_1,\Gamma_2$ are $(-4)$-curves with $\Gamma_1\Gamma_2=4$.
\item [(4)]  $F_{\infty}$ is an irreducible rational curve with $3$ nodes.
\end{itemize}
A direct computation gives the relative invariants \[K_{f_0}^2 = 9, \chi_{f_0} = 3, e_{f_0} = 27.\]
To obtain a semi-stable reduction, apply the degree-two base change $\tau: \mathbb{P}^1\to \mathbb{P}^1$ totally ramified over $t=\pm 2$. The fiber $F_2$ becomes a smooth hyperelliptic curve, and we obtain a semi-stable fibration $f: X \to \mathbb{P}^1 $ with exactly $5$ singular fibers, and \[K_{f}^2 = 10, \chi_{f} = 3, e_{f} = 26.\]

The semi-stable fibration has one fiber with twelve nodes and four
fibers with $4,4,3,3$ nodes. Indeed, the six nodes of $F_{-2}$ each
acquire a $(-2)$-bridge under the double cover, while the fibers
over $-1$ and $\infty$ each lift to two unramified copies.

For this fibration, we have
$(g,K_f^2,e_f)=(3,10,26)$, with five singular fibers carrying
$(12,4,4,3,3)$ nodes.  Total ramification of degree $d$ at $k$
selected singular fibers and $\ell$ smooth fibers gives, with
$m=k+\ell$ and $N$ selected nodes,
\[
 A=(8m-6)d-8m,\qquad E=(4m+18)d-4m,\qquad
 B=E-N(d-1/d).
\]
The QS certificate is $D=A-E+N(d-6/5)$.  The first three choices below will be used in the proof of
Theorem~\ref{thm:full-rational-interval}; the fourth gives a further
interpolation interval.

\begin{proposition}\label{prop:g3-four-bridges}
For each row of Table~\ref{tab:g3-four-bridges}, every rational
number strictly between the indicated endpoint slopes, as well as
$r_0$ itself, occurs as the auxiliary orbifold Chern slope
of a surface with $s_2<0$ and big cotangent bundle.
\end{proposition}
\begin{table}[htbp]\centering\papertablestyle
\caption{Four genus-three interpolation certificates.
Selected fiber sizes are their node counts; $D>0$ certifies QS bigness.}
\label{tab:g3-four-bridges}
\begin{tabular}{rccrccc}\toprule
$d$ & selected nodes & $\ell$ & $N$ & $D$ & $r_0=A/B$ & $r_1=(2E-A)/(B+E-A)$\\\midrule
$19$ & $(3,3)$ & $3$ & $6$ & $54/5$ & $1919/1863$ & $2527/2167$\\
$21$ & $(12,3)$ & $1$ & $15$ & $33$ & $1239/1063$ & $3087/1987$\\
$9$ & $(12,4,3)$ & $0$ & $19$ & $141/5$ & $621/401$ & $1701/941$\\
$4$ & $(12,4,4,3,3)$ & $1$ & $26$ & $244/5$ & $80/31$ & $112/47$\\\bottomrule
\end{tabular}
\end{table}
\begin{proof}
The four values of $m=k+\ell$ are $5,3,3,6$.  Each admits a
connected degree-$d$ cover totally ramified at exactly the specified
points: for even $m$, use pairs of inverse $d$-cycles; for odd $m$
(the listed degrees are odd), use inverse pairs and the triple
$\sigma,\sigma,\sigma^{-2}$.  The generated action is transitive,
so the Riemann existence theorem applies.

Direct substitution gives $A>0$, $A<E$, $B>0$ and $D>0$
in every row; each base curve has genus at least two.  By
Proposition~\ref{prop:aow-certificate}, the smooth pullback has
$s_2<0$ and big cotangent bundle.  Proposition~\ref{prop:aow-interp}
then realizes all the stated rational values between $r_0$ and
$r_1$.  Only $r_1$ is excluded, because $s_2<0$ is strict.
\end{proof}

\subsection{Genus-two mixed ramification}\label{sec:genus2-mixed}
\label{subsec:six-fibers}
We use the six-fiber construction in
\cite[Example~5.2]{LTYZ14}, following the original set-up of
S.-L. Tan. Let $pr_i: \mathbb{P}^1 \times \mathbb{P}^1 \to \mathbb{P}^1$ be the projection onto the $i$-th factor. Take an irreducible smooth curve $B_0$ of type $(2,1)$. There are two fibers $\Gamma_1^0, \Gamma_2^0$ of $pr_2$ that are tangent to $B_0$. 

Take two other general fibers $\Gamma_3^0, \Gamma_4^0$ of $pr_2$. We construct a double cover $\pi_0: S_0\to \mathbb{P}^1 \times \mathbb{P}^1$ with $B_{\pi_0} = \sum_{i=1}^4 \Gamma_i^0$ as the branch locus. 

Let $B_\pi=\pi_0^*B_0$, with $\pi_0^*\Gamma_i^0=2\Gamma_i$. 
The first covering surface $S_0$ is isomorphic to $E\times\mathbb P^1$,
 where $E$ is the elliptic double cover obtained by pulling back a general fiber of $pr_1$. Let $H$ be a general elliptic fiber of its projection to $\mathbb P^1$, and choose $L=H+\Gamma_3$. Then $B_\pi\sim2L$, so this chosen square root defines the second double cover $\pi:S_1\to E\times\mathbb P^1$. Let $S$ be the minimal resolution of $S_1$. For $j=1,\ldots,6$, put $q_j=E_j\cap\Gamma_3$, so $L|_{E_j}\simeq\mathcal O_{E_j}(q_j)$. Put $p_i=E_i\cap\Gamma_i$ for $i=1,2$, $p_j=E_j\cap\Gamma_3$ for $j=3,4$, and $p_j=E_j\cap\Gamma_4$ for $j=5,6$, with the fibers grouped as below.

This induces a semi-stable genus-two fibration $f:S\to\mathbb P^1$ with six singular fibers, exactly four of which have non-compact Jacobians. The verification is as follows:

\begin{enumerate}
    \item Let $E_1, E_2$ be two fibers of $pr_1$ whose images pass through the tangent points of $B_0$ with $\Gamma_1^0, \Gamma_2^0$ respectively. The restrictions satisfy $B_\pi|_{E_i}=2p_i$ and $L|_{E_i}\simeq\mathcal O_{E_i}(q_i)$, hence $2p_i\sim2q_i$. The distinct ramification points $p_i,q_i$ differ by a nonzero two-torsion point of $E$. Consequently the normalization of $\pi^{-1}E_i$ is a connected unramified double cover, so $\pi^{-1}E_i$ is irreducible. Thus the corresponding fiber of $f$ is $\widetilde{F_1} = C_1 + C_2$ (similarly $\widetilde{F_2} = C_1 + C_2$), where $C_1 C_2 = 2$, $C_1$ is a smooth elliptic curve, and $C_2$ is a $(-2)$-curve.
    
    \item Let $E_3,E_4$ be the two elliptic fibers whose images pass through the two intersection points of $B_0$ with $\Gamma_3^0$. Since $p_3=q_3$ (similarly $p_4=q_4$), the inverse image
$\pi^{-1}E_3$ (respectively $\pi^{-1}E_4$) is
\emph{reducible}.  The corresponding fiber $F_3$ (respectively
$F_4$) consists of two smooth elliptic curves meeting transversely
at one point.  It has one node and compact Jacobian; cf.
\cite[Example~5.2]{LTYZ14}.
    
    \item Let $E_5,E_6$ be the two elliptic fibers whose images pass through the two intersection points of $B_0$ with $\Gamma_4^0$.

 The inverse images of $E_5,E_6$ are instead
\emph{irreducible}. 
Each of the corresponding fibers $F_5,F_6$ is an
irreducible nodal curve with elliptic normalization, and hence has one
node and a non-compact Jacobian
\cite[Example~5.2]{LTYZ14}.
\end{enumerate}

Consequently, the six node counts are $(2,2,1,1,1,1)$;
precisely $F_1,F_2,F_5,F_6$ have non-compact Jacobians.
We shall use these fiber types in the mixed base changes below.
Furthermore, we have $K_f^2 = 4$, $\chi_f = 1$, $e_f = 8$, $p_g(S) = 0$, and $q(S) = 1$.

The second base change used in Proposition~\ref{prop:aow-interp} was
branched only above \emph{smooth} fibers.  Allowing controlled
ramification over chosen \emph{nodal} fibers gives a considerably
stronger result.  Unlike that smooth-fiber interpolation, this
construction changes the types of the quotient singularities from
$A_{d-1}$ to $A_{de-1}$ at the additionally ramified points.  We
continue to contract only the exceptional chains created by the base
changes, leaving the $(-2)$-components of the designated initial
fibration untouched.

\subsubsection{The local formula}
\begin{lemma}
\label{lem:mixed-g2}
Let $Y\to C$ be a semi-stable genus-two fibration with $g(C)\geq2$,
obtained from a semi-stable fibration over $\mathbb P^1$ by a first
base change of degree $d\geq2$ which is totally ramified over the selected nodal fibers and unramified over every other singular-fiber point. Additional ramification over smooth-fiber points is allowed. Assume that the selected fibers contain $N\geq1$ nodes in total
before the first base change, and set
\[
 A=c_1^2(Y),\quad E=c_2(Y),\quad S=A-E,\qquad
 B=E-N\left(d-\frac1d\right).
\]
Choose $q$ of the already ramified fiber points on $C$ such that the
corresponding fibers of the original fibration contain a total of $T$
nodes (so $T\leq N$).  Take a
connected cover $C'\to C$ of degree $n$ whose local monodromy at
each chosen point has cycle type $e^k1^{\,n-ke}$, with
$e\geq2$, $k\geq0$, and $ke\leq n$, and which is unramified
elsewhere.  Assume that the local monodromies are paired with their
inverses, as will be the case for $q=2$ or $4$ below.  Let $Y'$ be the
smooth semi-stable pullback, and contract all the chains created by
the first and second changes to obtain $Z'_{\mathrm{aux}}$.
Then, with $R=qk(e-1)$,
\begin{align}
 c_1^2(Y')&=nA+4R, & c_2(Y')&=nE+2R,\label{eq:g2-mixed-AE}\\
 s_2(Y')&=nS+2R, &
 \frac{c_2(\mathcal X'_{\mathrm{aux}})}n
 &= B+\frac{2R}{n}-\frac{Tk}{dn}
                 \left(e-\frac1e\right).
 \label{eq:g2-mixed-B}
\end{align}
The auxiliary quotient singularities consist of
$Tk$ points of type $A_{de-1}$ and
$nN-keT$ points of type $A_{d-1}$; in particular, the number
of singularities is $kT+nN-keT$.
  If
\begin{equation*}
 D_0:=S+N\left(d-\frac65\right)\geq0,
\end{equation*}
then the Asega--De Oliveira--Weiss sufficient inequality for $Y'$ is
strictly satisfied, since
\begin{equation}\label{eq:g2-mixed-QS}
 s_2(Y')+6\sum_{z\in\mathrm{Sing}\,Z'_{\mathrm{aux}}}
 h^1_\Omega(z)
 >nD_0+2R+\frac65\,Tk(e-1)\geq0.
\end{equation}
\end{lemma}
\begin{proof}
The genus-two case of the Riemann--Hurwitz and semi-stable
base-change formulas gives~\eqref{eq:g2-mixed-AE}.  At an original
node, the equation $xy=t$ becomes $xy=u^d$ after the first covering.
Above an unramified point of the second covering this still produces
an $A_{d-1}$ chain; above a point with ramification index $e$ the
composite local equation is $xy=v^{de}$, producing an $A_{de-1}$
chain.  Each of the $T$ selected original nodes therefore gives $k$
chains of type $A_{de-1}$ and $n-ke$ of type $A_{d-1}$; the
remaining $N-T$ original nodes give $n$ chains of type $A_{d-1}$.
These chains can be contracted while retaining the components of the
designated initial fibration.

For a contracted $A_{r-1}$ chain, the reduction in the orbifold
Euler number relative to the smooth resolution is $r-1/r$.
Consequently,
\begin{align*}
 c_2(\mathcal X'_{\mathrm{aux}})
 &=c_2(Y')-(nN-keT)\left(d-\frac1d\right)
       -kT\left(de-\frac1{de}\right)\\
 &=nB+2R-\frac{Tk}{d}\left(e-\frac1e\right),
\end{align*}
which proves~\eqref{eq:g2-mixed-B}.
Finally use $6h^1_\Omega(A_{r-1})>r-6/5$ for every $r\geq2$.
Replacing $e$ ordinary $A_{d-1}$ points by one $A_{de-1}$ point
increases this rational lower bound by
$\frac65(e-1)$; there are $Tk$ replacements.
Together with $s_2(Y')=nS+2R$ this gives
\eqref{eq:g2-mixed-QS}.
\end{proof}

\begin{remark}
\label{rem:g2-mono}
For $q=2$ or $4$ and every rational
$\theta\in[0,1]$, choose $e\geq2$, write $\theta=a/b$ with
$0\leq a\leq b$, and set $n=eb$, $k=a$.  Let $\tau\in S_n$ have
cycle type $e^k1^{\,n-ke}$.  Assign the peripheral loops around the
chosen points the permutations
$\tau,\tau^{-1}$, or
$\tau,\tau^{-1},\tau,\tau^{-1}$, respectively.  As $g(C)\geq2$,
we may send one handle generator to an $n$-cycle and all remaining
handle generators to the identity.  The surface-group relation is
satisfied and the monodromy action is transitive.  The Riemann
existence theorem supplies a connected cover with exactly these
ramification data and no other branch points.  The proportion of
sheets in the $e$-cycles is precisely
$\theta=ke/n$.  For $k=0$, the cover is connected and \emph{\'etale}.
\end{remark}

\subsubsection{The intervals around two}
For the fibration with six singular fibers above, the relative invariants are
$K_f^2=4$ and $e_f=8$.  The singular fibers have, in order,
\[
                 (2,2,1,1,1,1)
\]
nodes. The first two fibers each have an elliptic component and an
old $(-2)$-component meeting at two nodes; $F_3,F_4$ are reducible,
with two smooth elliptic components meeting at one point, whereas
$F_5,F_6$ are irreducible one-node fibers.  The latter four fibers
each have one node.  We do not contract these old $(-2)$-components in the
auxiliary construction.

\begin{theorem}
\label{thm:g2-filled}
For every rational
\[
          r\in\left[\frac53,\frac{408}{169}\right],
\]
there is a projective normal auxiliary surface $Z_{\mathrm{aux}}$
with only Du Val singularities such that its smooth minimal
resolution $Y$ is a minimal surface of general type satisfying
$s_2(Y)<0$ and having big cotangent bundle, with
$c_1^2(\mathcal X_{\mathrm{aux}})/c_2(\mathcal X_{\mathrm{aux}})=r$.
Only chains introduced by the base changes are contracted.
\end{theorem}
\begin{proof}
First take the degree-three cover totally ramified above
$F_1,F_2,F_5,F_6$, the four fibers with non-compact Jacobians;
they contain $N=2+2+1+1=6$ original nodes.
The resulting base curve has genus two, and the smooth pullback has
$(A,E,S,B)=(20,28,-8,12)$ and $D_0=14/5>0$.
Choose the second cover of Lemma~\ref{lem:mixed-g2} with the prescribed
local monodromy at all four corresponding points, so $q=4$ and $T=6$.
 Put
$\theta=ke/n\in[0,1]\cap\mathbb Q$ and
$u=\theta(1-1/e)$, $v=\theta(1-1/e^2)$.  Equations
\eqref{eq:g2-mixed-AE}--\eqref{eq:g2-mixed-B} yield
\begin{equation*}
 r_-(\theta,e)=\frac{20+16u}{12+8u-2v}.
\end{equation*}
The smooth Segre number and QS lower certificate are
\begin{equation}\label{eq:g2-minus-QS}
 s_2(Y')=-8n+8k(e-1)<0,\qquad
 nD_0+2R+\frac65Tk(e-1)
 =\frac{14n+76k(e-1)}5>0.
\end{equation}
For each fixed $e\geq2$, the function $\theta\mapsto
r_-(\theta,e)$ is a strictly increasing fractional-linear map with
rational coefficients.  Its endpoint values are
\[
 r_-(0,e)=\frac53,\qquad
 r_-(1,e)=\frac{12(3e)-16}
 {6(3e)-8+\frac6{3e}}<2,
 \qquad \lim_{e\to\infty}r_-(1,e)=2.
\]
Thus for every rational $r\in[5/3,2)$ one can choose $e$ such
that $r\leq r_-(1,e)$, solve uniquely for a rational
$\theta\in[0,1]$, and realize it using
Remark~\ref{rem:g2-mono}.  Every such $Y'$ satisfies the required
Segre and QS conditions by~\eqref{eq:g2-minus-QS}.

For slopes at least two, start instead with the degree-two cover
ramified over \emph{all six} singular fibers.  The base curve again
has genus two, and
$(A,E,S,B)=(16,20,-4,8)$ with $D_0=12/5>0$.
Now branch the second covering only over the first two fibers,
which each contain two original nodes.  Hence $q=2$ and $T=4$.
Using the same $\theta,u,v$,
\begin{equation*}
 r_+(\theta,e)=\frac{16+8u}{8+4u-2v},\qquad
 s_2(Y')=-4n+4k(e-1)<0,
\end{equation*}
and the QS lower certificate is
\begin{equation*}
 nD_0+2R+\frac65Tk(e-1)
 =\frac{12n+44k(e-1)}5>0.
\end{equation*}
Fixing $e=6$ gives the particularly simple interpolation
\begin{equation*}
 \displaystyle
 r_+(\theta,6)=\frac{288+120\theta}{144+25\theta},
 \qquad 0\leq\theta\leq1.
\end{equation*}
This is strictly increasing and carries $[0,1]\cap\mathbb Q$
bijectively onto
$[2,408/169]\cap\mathbb Q$.
The endpoint $\theta=1$ is realized with $n=e=6$ and $k=1$;
it gives
\[
\begin{aligned}
 c_1^2(Y')&=136,& c_2(Y')&=140,\\
 s_2(Y')&=-4,& c_2(\mathcal X'_{\mathrm{aux}})&=\frac{169}{3}.
\end{aligned}
\]
Here $Z'_{\mathrm{aux}}$ has precisely four $A_{11}$ points and
24 $A_1$ points, obtained by contracting only the base-change chains;
the four $A_{11}$ points arise over the two chosen double-node
fibers.  Thus the endpoint slope is $136/(169/3)=408/169$.
The two constructions join at slope~$2$, establishing the theorem.
\end{proof}

\subsection{Proof of Theorem~\ref{thm:full-rational-interval}}
We conclude the proof of Theorem~\ref{thm:full-rational-interval}.

\begin{proof}
Proposition~\ref{prop:tu-near-one} realizes all rational slopes in
$(1,33/32)$.  The first three rows of
Table~\ref{tab:g3-four-bridges} extend this interval beyond $5/3$:
indeed,
\[
 \frac{1919}{1863}<\frac{33}{32},\qquad
 \frac{1239}{1063}<\frac{2527}{2167},\qquad
 \frac{621}{401}<\frac{3087}{1987},\qquad
 \frac53<\frac{1701}{941}.
\]
Theorem~\ref{thm:g2-filled} then gives
$[5/3,408/169]\cap\mathbb Q$. Finally,
Corollary~\ref{cor:sharp-three-interval} gives
$(99/43,3]\cap\mathbb Q$. These two intervals overlap because
\[
 \frac{99}{43}<\frac{408}{169}.
\]
All the overlaps are strict. Hence the intervals cover
$(1,3]\cap\mathbb Q$, including all the intermediate rational endpoints.
\end{proof}

\section{Further canonical models and constructions}\label{sec:further}
We first determine the optimal canonical-model orbifold Chern slope
for the genus-three pencil used in Section~\ref{sec:aow-density}.
We then give asymptotic Segre inequalities and further examples from
semi-stable fibrations. The corresponding interpolation under the
Roulleau--Rousseau criterion was recorded in
Corollary~\ref{prop:interp}.
\subsection{An optimal canonical-model slope for a genus-three pencil}\label{subsec:g3-canonical}
In~\cite[Section~5]{RR14}, Roulleau and Rousseau considered the
geography of smooth surfaces $Y$ with $s_2(Y)<0$ and big cotangent
bundle.  Under their orbifold Segre-sum hypothesis
$s_2(Y)+s_2(\mathcal X)>0$, they obtained the numerical restriction
\[
     \frac{c_1^2(Y)}{c_2(Y)}>\frac35.
\]
They also observed that examples near this lower bound can give
orbifold Chern slopes approaching the Bogomolov--Miyaoka--Yau
bound~$3$; see~\cite[Section~5]{RR14}.

\begin{proposition}\label{thm:near-three}
For the genus-three semi-stable pencil of Subsection~\ref{sub2}, a
finite connected base change produces a smooth minimal surface $Y$
of general type with $s_2(Y)<0$ and big cotangent bundle such that
its canonical-model orbifold $\mathcal X_{\rm can}$ has slope
\[
 \frac{c_1^2(\mathcal X_{\rm can})}
 {c_2(\mathcal X_{\rm can})}
 =\frac{5970}{1991}\approx 2.99849.
\]
Moreover, this value is the largest canonical-model orbifold Chern
slope among finite connected base changes of this fixed pencil
whose new base has genus at least two and whose canonical orbifold
has positive second Chern number.
\end{proposition}

\begin{lemma}\label{lem:oldbridges}
For the genus-three pencil of Subsection~\ref{sub2}, the semi-stable
reduction by the initial double cover has a singular fiber with $12$
nodes, including precisely six exceptional rational $(-2)$-bridges.
The other four singular fibers have $4,4,3,3$ nodes, respectively.
Take a connected degree-$d$ cover $C'\to\mathbb P^1$ with $g(C')\geq2$. Let $e_1,\ldots,e_r$ be its ramification indices above this $12$-node fiber, including the indices $1$ at unramified points; thus $\sum_j e_j=d$. If the auxiliary quotient contracts the new chains and retains the six old bridges, then
\[
 B_{\rm can}:=c_2(\mathcal X_{\rm can})
 =B_{\rm aux}-9\sum_{j=1}^r\frac1{e_j}.
\]
In particular, the correction is $9/d$ for total ramification and $9d$ when the cover is unramified over this fiber.
Contracting all these chains gives the canonical model of the
semi-stable pullback.
\end{lemma}
\begin{proof}
The original fiber $F_{-2}$ is a union of four rational components
meeting pairwise once, hence has six nodes.  The initial double base
change is ramified at $t=-2$.  At each original node its local model
becomes $xy=u^2$ and its resolution inserts one $(-2)$-bridge,
replacing the node by two nodes.  Thus this fiber has $12$ nodes,
including six old bridges.  The fibers over $t=-1,\infty$ lift without
ramification to two copies apiece, contributing $4,4,3,3$ nodes.

First suppose that the degree-$d$ cover is totally ramified at the special fiber. Each old bridge
and the two adjacent new chains combine into the resolution chain of
$xy=v^{2d}$, an $A_{2d-1}$ singularity.  The orbifold-Euler correction
from these six singularities is
\[
6\left(2d-\frac{1}{2d}\right)=12d-\frac3d,
\]
whereas contracting only the twelve newly created chains at its nodes
would contribute $12(d-1/d)=12d-12/d$.  Their difference is $9/d$.
If the special fiber is unramified, its six old $A_1$ bridges occur in
$d$ copies; each contributes $2-1/2$, and the difference is $9d$.
For arbitrary ramification, apply the same local computation with $d$ replaced by each local index $e_j$ and sum the corrections $9/e_j$.

The pullback $Y$ is minimal of general type, so $K_Y$ is nef
and big.  If an irreducible curve $D$ satisfies $K_Y\cdot D=0$,
then the Hodge index theorem gives $D^2<0$.  Adjunction gives
$2p_a(D)-2=D^2$, hence $D$ is a smooth rational $(-2)$-curve.
No such rational curve is horizontal, since the base has genus at
least two.  Thus $D$ lies in a singular fiber.  In a relatively
minimal semi-stable fibration, a smooth rational component meeting the other
components at $\nu$ nodes satisfies $D^2=-\nu$ and
$K_Y\cdot D=\nu-2$.  The original smooth rational stable components have $\nu\geq3$; the old and new bridges are exactly the components
with $\nu=2$.
  The fiber description in Subsection~\ref{sub2}
therefore exhausts the $K_Y$-trivial curves.  Their contraction
is the canonical model.
\end{proof}

\begin{proof}[Proof of Proposition~\ref{thm:near-three}]
Label the five singular fibers by $F_0,\ldots,F_4$, where $F_0$
contains the six old $(-2)$-bridges and has $12$ nodes, $F_1,F_2$
have four nodes each, and $F_3,F_4$ have three nodes each.  Choose
a coordinate in which their images $p_0,\ldots,p_4$ are finite.
The cyclic cover of degree~$30$ given by the normalization of
\[
 z^{30}=(t-p_0)^6(t-p_1)^5(t-p_2)^{25}
                   (t-p_3)^6(t-p_4)^{18}
\]
is connected, since $\gcd(30,6,5,25,6,18)=1$.  The exponents sum
to~$60$, so there is no further branching at infinity.  Above
$F_0,F_3,F_4$ each ramification cycle has length~$5$, and above
$F_1,F_2$ each has length~$6$.  The total ramification is
\[
 R=3\cdot6(5-1)+2\cdot5(6-1)=122,
\qquad 2g(C')-2=-60+122=62.
\]
Hence $g(C')=32$.  The standard semi-stable base-change formulas for
the relative invariants $(K_f^2,e_f)=(10,26)$ give
\[
 \begin{aligned}
 A:=c_1^2(Y)&=8R-6\cdot30=796,\\
 E:=c_2(Y)&=4R+18\cdot30=1028,\qquad s_2(Y)=-232.
 \end{aligned}
\]
At a node with ramification index $e$, the exceptional chain is of
type $A_{e-1}$ and contributes $e-1/e$ to the orbifold Euler
correction.  At $F_0$, however, the two new chains adjoining each
old bridge combine with it to form an $A_{2e-1}$ chain.  Thus the
canonical-model correction from $F_0$ is
$6\cdot6(10-1/10)$; the auxiliary correction would instead come
from twelve independent $A_4$ chains over each of the six points
above $F_0$.  Including the other four fibers gives
\begin{align*}
 B:=c_2(\mathcal X_{\rm can})
 &=1028-36\left(10-\frac1{10}\right)
       -40\left(6-\frac16\right)\\
 &\hspace{2.2em}-36\left(5-\frac15\right)
   =\frac{3982}{15}.
\end{align*}
Consequently
\[
 \frac{c_1^2(\mathcal X_{\rm can})}{c_2(\mathcal X_{\rm can})}
 =\frac{796}{3982/15}=\frac{5970}{1991},\qquad
 s_2(Y)+s_2(\mathcal X_{\rm can})
 =-232+796-\frac{3982}{15}=\frac{4478}{15}>0.
\]
Theorem~\ref{RR1} proves bigness of $\Omega_Y^1$.

We prove optimality for every finite connected base change of this
fixed pencil, allowing arbitrary ramification above the five singular
fibers and additional branching above smooth fibers.  Let its degree
be~$n$, and write $e_{ij}\geq1$ for the cycle lengths over the point
corresponding to $F_i$; thus $\sum_j e_{ij}=n$ for every~$i$.
Write $R_{\rm sm}\geq0$ for the total ramification over smooth-fiber
points and put
\[
 q=\sum_{i=0}^4\#\{j\},\qquad
 R=5n-q+R_{\rm sm},\qquad
 A=8R-6n,\qquad E=4R+18n.
\]
For the canonical model, the total orbifold correction from the
special fiber is
$6(2e-1/(2e))=12(e-1/e)+9/e$ for each cycle of
length~$e$.  Applying this together with the ordinary nodal
corrections gives
\begin{equation}\label{eq:g3-arbitrary-canonical}
 \begin{aligned}
 B:=c_2(\mathcal X_{\rm can})
 &=E-\sum_{i=0}^4\sum_j\nu_i
       \left(e_{ij}-\frac1{e_{ij}}\right)
       -9\sum_j\frac1{e_{0j}},\\
 &\hspace{2em}(\nu_0,\ldots,\nu_4)=(12,4,4,3,3).
 \end{aligned}
\end{equation}
The analysis of the $K_Y$-trivial curves in
Lemma~\ref{lem:oldbridges} applies cycle by cycle as soon as the
new base has genus at least two.  In particular,
\eqref{eq:g3-arbitrary-canonical} gives the second Chern number of
the \emph{canonical-model} orbifold, with the old bridges included
in the correction.

Substituting the preceding formulas and using
$\sum_j e_{ij}=n$ yields the exact identity
\begin{align*}
 5970B-1991A
 &=7952R_{\rm sm}
   +\sum_{i\in\{0,3,4\}}\sum_j
       \frac{(e_{ij}-5)(874e_{ij}-3582)}{e_{ij}}\\
 &\quad+\sum_{i\in\{1,2\}}\sum_j
       \frac{(e_{ij}-6)(662e_{ij}-3980)}{e_{ij}}.
\end{align*}
Each summand is nonnegative for every positive integer~$e_{ij}$:
in the first sum the two factors are both negative for $e\leq4$, their product vanishes at $e=5$, and they are both positive for $e\geq6$. In the second sum they are both negative for $e\leq5$, their product vanishes at $e=6$, and they are both positive for $e\geq7$.

Thus $5970B\geq1991A$.  Whenever $B>0$, this implies
$A/B\leq5970/1991$.  Equality requires $R_{\rm sm}=0$,
ramification index~$5$ at $F_0,F_3,F_4$, and index~$6$ at
$F_1,F_2$.  The cyclic degree-$30$ cover constructed above attains
these conditions and has $B>0$, proving the claim.
\end{proof}

The slopes in Corollary~\ref{cor:tu-strict-near-three} exceed this
fixed-pencil maximum for every $m\geq34$. Indeed,
\[
 3-\frac{2}{39m+2}>\frac{5970}{1991}
 \quad\Longleftrightarrow\quad 117m>3976.
\]

Take a connected degree-$d$ base change totally ramified at all five
singular fibers, where $d\geq3$ is odd.
The smooth semi-stable pullback
satisfies
\[
c_1^2(\Tilde{X})=34d-40,\quad c_2(\Tilde{X})=38d-20,
\quad s_2(\Tilde{X})=-4d-20<0.
\]
The auxiliary and canonical-model orbifold Euler numbers are different:
\[
c_2(\mathcal X_{\rm aux})=12d-20+\frac{26}{d},\qquad
c_2(\mathcal X_{\rm can})=12d-20+\frac{17}{d}.
\]
Indeed, the initial semi-stable double cover creates six older $(-2)$-bridges,
and the correction is $9/d$; see Lemma~\ref{lem:oldbridges}.
In particular,
\[
s_2(\Tilde X)+s_2(\mathcal X_{\rm can})
 =18d-40-\frac{17}{d}>0,
\quad
\frac{c_1^2(\mathcal X_{\rm can})}{c_2(\mathcal X_{\rm can})}
 =\frac{34d-40}{12d-20+17/d}.
\]
Theorem~\ref{RR1} therefore gives
bigness of $\Omega^1_{\Tilde X}$ with $s_2(\Tilde X)<0$.

\subsection{Asymptotic Segre inequalities under base changes}\label{subsec:asymptotic-segre}
Let $f:X\to\mathbb P^1$ be a non-isotrivial semi-stable fibration, and let $\mathcal R\subset\mathbb P^1$ be its set of singular values. We take a finite base change $\pi: \Tilde{C}\to \mathbb{P}^1$ of degree $d$ totally ramified over $\{p_1=f(F_1), p_2=f(F_2), \cdots, p_k=f(F_k), k\leq s\}\subset \mathcal{R}$ and over $\ell$ additional points of $\mathbb P^1\setminus\mathcal R$, and unramified elsewhere. This yields a non-isotrivial semi-stable fibration $\Tilde{f}: \Tilde{X}\to \Tilde{C}$. Denote the genus of $\Tilde{C}$ by $\Tilde{b}$. By Lemma~\ref{bc}, contracting the new $(-2)$-curves yields an auxiliary normal surface $\Hat X$ with $A_{d-1}$ singularities, which is of general type provided the smooth pullback is of general type (see the following commutative diagram).
\begin{center}
\begin{tikzcd}
\Hat{X} \arrow[r,  leftarrow] & \Tilde{X} \arrow[r] \arrow[d, "\Tilde{f}"]& X \arrow[d,  "f"] \\
    & \Tilde{C} \arrow[r, "\pi"]
      & \mathbb{P}^1

\end{tikzcd}
\end{center}
 We write $\mathcal X$ for the orbifold associated with the
quotient surface $\Hat X$.  The relevant Chern and Segre-number
formulas are collected in Lemma~\ref{S2} and in
\eqref{eq:AE-general}--\eqref{eq:aux-general}.

The Segre-sum test uses $s_2(\Tilde{X})<0$ together with
$s_2(\Tilde{X})+s_2(\mathcal X_{\rm aux})>0$.  The second
inequality is sufficient, but not necessary, for the stronger
QS-bigness certificate of Proposition~\ref{prop:aow-certificate}.

\begin{proposition}\label{prop:asymptotic-segre}
Let
$f\colon X\longrightarrow \mathbb P^1$
be a non-isotrivial semi-stable fibration of genus $g\geq 2$ with
$s\geq 5$ singular fibers. Let
$\pi_d\colon \widetilde C_d\longrightarrow \mathbb P^1
$
be a sequence of finite base changes of degree $d$,
totally ramified over all singular values of $f$ and unramified
elsewhere. Let $\widetilde X_d$ be the smooth semi-stable pullback,
and let $\mathcal X_d$ be the auxiliary orbifold obtained by
contracting the new exceptional chains. Set
\[
\rho_f:=\frac{(g-1)(s-2)}{\chi_f}.
\]
Then
\begin{align}
s_2(\widetilde X_d)
&=
2d\chi_f\bigl(\lambda_f-6+\rho_f\bigr)
-2s(g-1),
\label{eq:s2-Xtilde-asymptotic}
\\
s_2(\widetilde X_d)+s_2(\mathcal X_d)
&=
d\chi_f\bigl(3\lambda_f-12+4\rho_f\bigr)
-4s(g-1)
+\frac{\chi_f(\lambda_f-12)}{d}.
\label{eq:s2-sum-asymptotic}
\end{align}
In particular,
\begin{align*}
\lim_{d\to\infty}
\frac{s_2(\widetilde X_d)}{2d\chi_f}
&=
\lambda_f-6+\rho_f,
\\
\lim_{d\to\infty}
\frac{s_2(\widetilde X_d)+s_2(\mathcal X_d)}
     {d\chi_f}
&=
3\lambda_f-12+4\rho_f.
\end{align*}

Consequently, for all sufficiently large admissible degrees $d$,
\[
s_2(\widetilde X_d)<0
\qquad\text{and}\qquad
s_2(\widetilde X_d)+s_2(\mathcal X_d)>0
\]
if and only if
\begin{equation}\label{eq:asymptotic-slope-window}
4-\frac{4}{3}\rho_f
<
\lambda_f
\leq
6-\rho_f.
\end{equation}
Equivalently,
\[
4-\frac{4(g-1)(s-2)}{3\chi_f}
<
\frac{K_f^2}{\chi_f}
\leq
6-\frac{(g-1)(s-2)}{\chi_f}.
\]
\end{proposition}

\begin{proof}
Here $k=s$, $\ell=0$, and $N_0=e_f$.  Lemma~\ref{S2} yields
\begin{align*}
 s_2(\widetilde X_d)
 &=d(K_f^2-e_f)+2(g-1)\bigl(d(s-2)-s\bigr),\\
 s_2(\widetilde X_d)+s_2(\mathcal X_d)
 &=d(2K_f^2-e_f)-\frac{e_f}{d}
   +4(g-1)\bigl(d(s-2)-s\bigr).
\end{align*}
The relative Noether formula $K_f^2+e_f=12\chi_f$, together
with $K_f^2=\lambda_f\chi_f$, gives
\[
 K_f^2-e_f=2\chi_f(\lambda_f-6),\qquad
 2K_f^2-e_f=\chi_f(3\lambda_f-12),\qquad
 e_f=\chi_f(12-\lambda_f).
\]
Substitution proves~\eqref{eq:s2-Xtilde-asymptotic} and
\eqref{eq:s2-sum-asymptotic}; dividing by $d$ gives the two limits.

It remains to check the boundary cases in the claimed equivalence.
If $\lambda_f-6+\rho_f=0$, then
$s_2(\widetilde X_d)=-2s(g-1)<0$.
If $3\lambda_f-12+4\rho_f=0$, the other expression becomes
\[
 s_2(\widetilde X_d)+s_2(\mathcal X_d)
 =-4s(g-1)+\frac{\chi_f(\lambda_f-12)}{d}<0
\]
for all sufficiently large $d$.  The signs of the leading
coefficients therefore give the weak upper and strict lower bounds
in~\eqref{eq:asymptotic-slope-window}, respectively.
\end{proof}

\begin{corollary}\label{cor:asymptotic-segre-canonical}
Under the assumptions of Proposition~\ref{prop:asymptotic-segre}, suppose that
\begin{equation*}
\lambda_f>\frac{12}{5}
\qquad\text{and}\qquad
\rho_f\leq 6-\lambda_f.
\end{equation*}
Then, for all sufficiently large $d$,
\[
s_2(\widetilde X_d)<0
\qquad\text{and}\qquad
s_2(\widetilde X_d)+s_2(\mathcal X_d)>0.
\]

Moreover, these conditions imply
\[
\lambda_f<4.
\]
\end{corollary}

\begin{proof}
The strict canonical class inequality of~\cite{Tan95} gives
$K_f^2<2(g-1)(s-2)$, or $\lambda_f<2\rho_f$.
Consequently,
$3\lambda_f-12+4\rho_f>5\lambda_f-12>0$,
while $\rho_f\leq6-\lambda_f$ implies
$\lambda_f-6+\rho_f\leq0$.
Proposition~\ref{prop:asymptotic-segre} gives the first assertion.
Finally, $\lambda_f/2<\rho_f\leq6-\lambda_f$
forces $\lambda_f<4$.
\end{proof}

\begin{remark}\label{lowslope}
\begin{enumerate}
\item These observations suggest starting the construction with a semi-stable fibration $f: X \to \mathbb{P}^1$ of genus $g$ with slope $\lambda_f<4$ and few singular fibers.

    \item Tan proved in \cite{Tan95} that a non-isotrivial semi-stable fibration $f:X\to\mathbb P^1$ of genus $g\geq2$ has at least five singular fibers, confirming a conjecture of Beauville. Moreover, if $\kappa(X)\geq0$ (resp. $\kappa(X)=2$), then $s\geq6$ (resp. $s\geq7$); see \cite{LTZ18,TTZ05}.
\end{enumerate}
\end{remark}

\subsection{Further genus-two examples}
We first use the five-fiber fibration of
Subsection~\ref{subsec:tu-five}.
For an odd-degree cover totally ramified over the five singular
fibers,
Lemma~\ref{S2}, with $(g,K_f^2,e_f,N_0,k,\ell)=(2,5,19,19,5,0)$,
gives
\begin{align*}
 c_1^2(\Tilde X)&=17d-20, & c_2(\Tilde X)&=25d-10,\\
 s_2(\Tilde X)&=-8d-10<0, &
 s_2(\Tilde X)+s_2(\mathcal X_{\rm aux})&=3d-20-\frac{19}{d}.
\end{align*}

For odd $d\geq9$, the Segre sum is positive.
The auxiliary Segre-sum inequality implies the QS condition by Remark~\ref{rem:aow-versus-rr}; hence Proposition~\ref{prop:aow-certificate} gives a surface of general type with negative second Segre number and big cotangent bundle.

For this five-fiber fibration, take a connected cover of odd degree $d$
totally ramified over the four singular values $p_0,p_1,p_2,p_3$ of
Lemma~\ref{lem:tu-canonical-chains} and one additional smooth-fiber
value, and unramified elsewhere. Then $N_0=18$, and Lemma~\ref{S2} gives
\[
 s_2(\Tilde X)+s_2(\mathcal X_{\rm aux})=2d-20-\frac{18}{d}.
\]

For odd $d\geq11$ this sum is positive.
Hence, $\Omega^1_{\Tilde{X}}$ is big.

\vspace{2mm}

For comparison, write $Y=\widetilde X$ for the smooth pullback.
For an odd-degree cover totally ramified over all five singular
fibers, the smooth and auxiliary orbifold Chern slopes are
\begin{equation}\label{eq:geography-Tu}
 \frac{c_1^2(Y)}{c_2(Y)}=\frac{17d-20}{25d-10},\qquad
 \frac{c_1^2(\mathcal X_{\rm aux})}{c_2(\mathcal X_{\rm aux})}
 =\frac{17d-20}{6d-10+19/d}.
\end{equation}
For example, $d=9$ and $d=11$ give auxiliary slopes approximately
$2.884$ and $2.893$, while the corresponding smooth Chern slopes
are approximately $0.619$ and $0.630$.
  If instead the four singular fibers with non-compact Jacobians
and one smooth fiber are chosen, the auxiliary denominator
in~\eqref{eq:geography-Tu} is $7d-10+18/d$.

We next use the six-fiber fibration of Subsection~\ref{subsec:six-fibers}.
For the unbounded-degree family we branch only over the four Jacobian-singular fibers.
Total ramification over all six has $s_2=4d-12$, so the negative-Segre condition fails for $d\geq3$; however, the exceptional degree-two cover branched over all six has $s_2=-4$ and is crucial in Subsection~\ref{sec:genus2-mixed}. When we totally ramify over $F_1,F_2,F_5,F_6$, the four
non-compact-Jacobian fibers, $N_0=6$ and the base degree can be any integer $d\geq3$; by Lemma~\ref{S2}, 

\[\begin{aligned}
&c_1^2(\Tilde{X})=12d-16, ~ c_2(\Tilde{X})=12d-8, ~ s_2(\Tilde{X})=-8<0;\\
&c_1^2(\mathcal{X})=12d-16, ~ c_2(\mathcal{X}_{\rm aux})=6d-8+\frac{6}{d}, ~ s_2(\mathcal{X}_{\rm aux})=6d-8-\frac{6}{d}>0~(d\geq4).
\end{aligned}\]

Hence the auxiliary Segre sum equals $6d-16-\frac{6}{d}$ and is strictly positive for $d\geq4$. At $d=3$ it vanishes, but the stronger AOW rational certificate is $-8+6(3-6/5)=14/5>0$, so Proposition~\ref{prop:aow-certificate} still proves that $\Tilde{X}$ has negative second Segre number and big cotangent bundle. It follows that

\[ \frac{c_1^2(\mathcal{X})}{c_2(\mathcal{X})}=\frac{12d-16}{6d-8+\frac{6}{d}},
\qquad\text{and hence~}~ \frac{c_1^2(\mathcal{X})}{c_2(\mathcal{X})}\to 2, \quad \text{as $d\to\infty$}.\]

For comparison, an odd-degree base change totally ramified over the five points $p_j=f(F_j)$, $j=1,\ldots,5$, has $N_0=7$ and auxiliary orbifold Chern slope
\[
 \frac{c_1^2(\mathcal X_{\rm aux})}{c_2(\mathcal X_{\rm aux})}
 =\frac{16d-20}{7d-10+7/d}.
\]
If instead the selected points are $p_j=f(F_j)$, $j=2,\ldots,6$, then $N_0=6$ and
\[
 \frac{c_1^2(\mathcal X_{\rm aux})}{c_2(\mathcal X_{\rm aux})}
 =\frac{16d-20}{8d-10+6/d}\longrightarrow2
 \qquad(d\to\infty).
\]
Both five-point comparison schemes have $s_2(\Tilde X)=2d-10$;
thus for odd $d\geq5$ they do not satisfy the negative-Segre condition.

\subsection{Higher-genus fibrations}\label{sub3}
We also obtain higher-genus fiber spaces with negative second Segre number and big cotangent bundle.

In \cite[Theorem 1.1(1)(a)]{Tu07}, Y.-P. Tu constructed a semi-stable fibration 
$f: X \to \mathbb{P}^1 $ of genus $g=2^n-1\,(n\geq2)$ with $6$ singular fibers, and \[K_{f}^2 = 6g-6, \chi_{f} = g+1, e_{f} = 6g+18.\]

Relabel the six singular fibers $F_i$ in nonincreasing order of their numbers of nodes. Let $N_0$ denote the total number of nodes in $F_1,\ldots,F_4$. Then $N_0\geq \frac{2}{3}e_f=4g+12$. 

A degree-two base change totally ramified over $p_i=f(F_i)$, $i=1,\ldots,4$, yields a semi-stable fibration $\Tilde{f}: \Tilde{X}\to \Tilde{C}$ over an elliptic curve $\Tilde{C}$. 
Since the pullback fibration is relatively minimal and every rational curve maps constantly to the elliptic base, $\Tilde X$ is minimal. The generically finite map to $X$, which is birational to a K3 surface, implies $\kappa(\Tilde X)\geq0$. The equality $K_{\Tilde X}^2=12g-12>0$ computed below then shows that $\Tilde X$ is of general type. By Lemma~\ref{bc}, the base change produces $N_0$ $(-2)$-curves. Contracting these $(-2)$-curves gives $\Hat{X}$ with $N_0$ singularities of type $A_1$. Hence
\[c_1^2(\Tilde{X})=2K^2_f=12g-12,\, c_2(\Tilde{X})=12g+36,\, s_2(\Tilde{X})=-48<0. \]

By Proposition~\ref{obd},
\[c_1^2(\mathcal{X})=2K^2_f=12g-12,\] \[ c_2(\mathcal{X})=12g+36-2N_0+\frac{N_0}{2}=12g+36-\frac{3}{2}N_0. \]
Since $4g+12\leq N_0\leq e_f=6g+18$, we have $0<3g+9\leq c_2(\mathcal X)\leq6g+18$.
It follows that $s_2(\mathcal{X})\geq 6g-30$ and \[\frac{c_1^2(\mathcal{X})}{c_2(\mathcal{X})}\geq \frac{12g-12}{6g+18}=\frac{2g-2}{g+3}\to 2.\]

Then \[s_2(\mathcal{X})+s_2(\Tilde{X})\geq 6g-30-48>0,\, \text{if}\,\, g>13.\]

The same auxiliary Segre-sum inequality implies the QS condition of Proposition~\ref{prop:aow-certificate}, so for $g=2^n-1>13$ (that is, $n\geq4$), $\Tilde{X}$ has negative second Segre number and big cotangent bundle.

For the genus-three specialization $g=3$ of the original six-fiber family, take a connected degree-$d$ cover totally ramified over all six singular-fiber points and unramified elsewhere, with $d\geq2$. Contracting the resulting new chains gives the auxiliary slope
\[
\frac{44d^2-48d}{16d^2-24d+36}.
\]
Among integers $d\geq2$, its unique maximum occurs at $d=10$ and equals
$980/349$. Indeed,
\[
 \frac{980}{349}-\frac{44d^2-48d}{16d^2-24d+36}
 =\frac{36(d-10)(9d-98)}{349(16d^2-24d+36)}\geq0,
\]
because the two numerator roots are $10$ and $98/9\in(10,11)$,
and the denominator is positive. 

\vspace{0.2cm}

\begin{samepage}
\medskip
\noindent{\bf Acknowledgements}
\phantomsection
\addcontentsline{toc}{section}{Acknowledgements}
The authors are supported by the NSF of China. X. L\"u is partially supported by the Scientific Research Innovation Capability Support Project for Young Faculty (No. SRICSPYF-ZY2025161), the Shanghai Pilot Program for Basic Research (No. TQ20240202), and the Science and Technology Commission of Shanghai Municipality (No. 22DZ2229014). W.-Y. Xu (corresponding author) is partially supported by the Science and Technology Innovation Plan of Shanghai (Grant No. 23JC1403200). The authors would like to thank Prof. Yifei Chen for useful conversations.

\vspace{0.2cm}
\noindent\textbf{AI declaration.}
The research questions were formulated, and the main mathematical proofs were developed and written by the authors in 2025. AI (GPT Astra) was used to assist with literature searches, English-language editing, and certain direct computations. The authors have independently verified the mathematical arguments and references and take full responsibility for the accuracy and content of the paper.

\end{samepage}

\phantomsection

\vspace{.1cm}
\phantomsection
\addcontentsline{toc}{section}{Contact}

\small
\noindent

{\it Address }: Xin L\"u: xlv@math.ecnu.edu.cn;
School of Mathematical Sciences, Ministry of Education Key Laboratory of Mathematics and Engineering Applications \& Shanghai Key Laboratory of PMMP,  East China Normal University, Shanghai 200241, China

{\it Address }: W.-Y. Xu: wanyuanxu@shnu.edu.cn;
Department of Mathematics, Shanghai Normal University, Shanghai 200234, P. R. of China

\end{document}